\documentclass{gtart}
\usepackage{amsmath, amssymb, graphicx, graphs, verbatim, epsf,amscd}

\newcommand\Span{\mathrm{Span}}
\newcommand\Pd{{\mathfrak {p}}}
\newcommand\ModFlow{{\widehat {\mathfrak {M}}}}

\newcommand\Mas{\mu}
\newcommand\Fa{\widehat F}
\newcommand\Sym{\mathrm{Sym}}
\newcommand\Aut{\mathrm{Aut}}
\newcommand\Page{P}
\newcommand\MCG{\mathrm{MCG}}
\newcommand\DiffP{\mathrm{Diff}^+}
\newcommand\Zmod[1]{\Z/{#1}\Z}
\newcommand\smargin[1]{\marginpar{\tiny{#1}}}
\newcommand{\alphas}{\mathbf{\alpha}}
\newcommand{\betas}{\mathbf{\beta}}
\newcommand{\gammas}{\mathbf{\gamma}}
\newcommand{\deltas}{\mathbf{\delta}}
\newcommand\TransInv{\mathfrak{T}}
\newcommand\TransInva{\widehat{\mathfrak{T}}}
\newcommand\LegInv{\mathfrak{L}}
\newcommand\LegInva{\widehat{\mathfrak L}}

\newcommand{\HFKa}{\widehat{\mathrm{HFK}}}
\newcommand\spinct{\mathbf t}

\newcommand\Field{\mathbb{F}}

\newcommand{\spinc}{\mathfrak{s}}
\newcommand\Ta{\mathbb{T}_{\alpha}}
\newcommand\Tb{\mathbb{T}_{\beta}}
\newcommand\Tc{\mathbb{T}_{\gamma}}
\newcommand\Td{\mathbb{T}_{\delta}}
\newcommand\CFa{\widehat{{\mathrm {CF}}}}
\newcommand\CFm{{{\mathrm {CF}}}^-}
\newcommand\HFa{\widehat{\mathrm {HF}}}
\newcommand\HFm{{{\mathrm {HF}}}^-}
\newcommand{\HFKm}{\mathrm{HFK}^-}
\newcommand{\CFKm}{\mathrm{CFK}^-}
\newcommand{\SpinC}{{\mbox{Spin}}^c}
\newcommand\RelSpinC{\underline{\SpinC}}

\newtheorem{thm}{Theorem}[section]
\newtheorem{cor}[thm]{Corollary}
\newtheorem{lem}[thm]{Lemma}
\newtheorem{prop}[thm]{Proposition}

\theoremstyle{definition}
\newtheorem{defn}[thm]{Definition}

\theoremstyle{remark}
\newtheorem{rem}[thm]{Remark}

\numberwithin{equation}{section}

\newcommand{\bfz}{{\mathbb {Z}}}

\newcommand{\x}{{\bf {x}}}

\newcommand{\y}{{\bf {y}}}
\newcommand{\s}{\mathbf s}

\newcommand{\Z}{\mathbb Z}

\newcommand{\Q}{\mathbb Q}

\newcommand{\bfr}{\mathbb R}

\DeclareMathOperator{\sln}{sl}
\DeclareMathOperator{\tb}{tb}
\DeclareMathOperator{\TB}{TB}

\DeclareMathOperator{\rot}{rot}

\DeclareMathOperator{\PD}{PD}

\begin{abstract}
  We study naturality properties of the transverse invariant in knot
  Floer homology under contact $(+1)$--surgery. This can be used as a
  calculational tool for the transverse invariant. As a consequence,
  we show that the Eliashberg--Chekanov twist knots $E_n$ are not
  transversely simple for $n$ odd and $n> 3$.
\end{abstract}


\begin{document}

\title[{Contact surgeries and the transverse invariant
in knot Floer homology}]
{Contact surgeries and the transverse invariant
in knot Floer homology}

\author[Peter Ozsv{\'a}th]{Peter Ozsv\'ath}
\address{Department of Mathematics, Columbia University,\\ 
New York 10027, USA}
\email{petero@math.columbia.edu}
\thanks{PSO was supported by NSF grant number DMS 0234311}

\author[Andr{\'a}s I. Stipsicz]{Andr{\'a}s I. Stipsicz}
\address{R\'enyi Institute of Mathematics, Budapest, Hungary and\\
  Department of Mathematics, Columbia University, New York 10027, USA}
\email{stipsicz@renyi.hu} 
\thanks{AS was supported by OTKA T49449, by CMI and by Marie Curie 
TOK project BudAlgGeo}
\maketitle

\section{Introduction}
Suppose that $T\subset (Y, \xi )$ is a null--homologous transverse
knot in the closed contact 3--manifold $(Y, \xi)$.  According
to~\cite{LOSS}, there is an invariant $\TransInv(T)$ of the transverse
isotopy class of $T$, taking values in the knot Floer homology group
$\HFKm(-Y,T)$ (introduced in~\cite{OSzknot, RasmussenThesis}). This
invariant is defined using open book decompositions and Heegaard Floer
homology.  For the definition, we approximate $T$ by a coherently
oriented Legendrian knot $L$ and find an appropriate open book
decomposition compatible with $(Y,\xi,L)$. In the Heegaard diagram
corresponding to this open book decomposition there is a distinguished
intersection point, giving a generator of the chain complex
$\CFKm(-Y,L)$ for knot Floer homology.  The element induces a homology
class $\LegInv(Y,\xi,L)\in\HFKm(-Y,L)$, called the \emph{Legendrian
  invariant} of $L$.  Since the Legendrian invariant remains unchanged
under negative stabilization, it can be viewed as an invariant
$\TransInv (Y,\xi,T)$ of the transverse knot $T$. This invariant turns out
to be an effective tool for studying Legendrian and transverse knots
in various contact 3--manifolds. The specialization $U=0$ turns $\HFKm
(-Y, L)$ into $\HFKa (-Y,L)$; the image of the invariant $\LegInv (Y,\xi,L)$
(and $\TransInv(Y,\xi,T)$)
under this reduction is denoted by $\LegInva (Y,\xi,L)$ (and $\TransInva(Y,\xi,T)$, resp.).

The motivation for the Legendrian and transverse invariants comes from
the construction in~\cite{OST}, which gives an invariant for
Legendrian and transverse knots in the standard contact three-sphere,
taking values in the combinatorial knot Floer homology of~\cite{MOS,
  MOST}; see also \cite{HKM2, OSzcont} for other constructions
and~\cite{NOT, VV} for related computations.

In this paper we show that the invariant $\LegInv$ (and hence
$\TransInv$) enjoys a simple transformation rule under the change
of the contact 3--manifold by contact $(+1)$--surgery. 

\begin{thm}\label{t:transz}
  Suppose that $L, S\subset (Y, \xi )$ are disjoint Legendrian knots
  in the contact 3--manifold $(Y, \xi )$, with $L$ 
  null--homologous and oriented.  Let $(Y_S, \xi _S)$ denote the
  contact 3--manifold we get by performing contact $(+1)$--surgery
  along $S$, while $L_S$ will denote the Legendrian knot $L$ viewed in
  $(Y_S, \xi _S)$. Suppose that $L_S$ is  null--homologous
  in $Y_S$.  The surgery gives rise to maps
 \[
F_{S,\spinc}\colon \HFKm(-Y, L)\to \HFKm (-Y_S, L_S),
\]
where $\spinc$ is a $\SpinC$ structure on the cobordism $W$ from $Y$
to $Y_S$.  There is a unique $\spinc$ for which 
\[
F_{S,\spinc}(\LegInv (Y,\xi, L))=\LegInv (Y_S, \xi_S, L_S)
\]
holds, and for all other $\SpinC$ structure $s$ the map $F_{S,s}$ is
trivial on $\LegInv (Y, \xi , L)$.  A similar identity holds for the
Legendrian invariant $\LegInva$ in $\HFKa$.
\end{thm}

This has the following immediate consequence for the transverse
invariant:

\begin{cor}
  Let $T\subset (Y, \xi )$ be a null--homologous transverse knot and
  $S\subset (Y, \xi )$ a Legendrian knot disjoint from $T$.  Let
  $(Y_S, \xi _S, T_S)$ denote the result of the contact
  $(+1)$--surgery along $S$, and suppose that $T_S$ (the knot $T$
  viewed in $(Y_S, \xi _S)$) is null--homologous in $Y_S$.  Then there
  is a unique $\SpinC$ structure $\spinc$ for which $F_{S, \spinc }
  (\TransInv (Y,\xi,T))=\TransInv (Y_S,\xi_S,T_S)$ and for all other $\SpinC$
  structure $s$ the map $F_{S,s}$ is trivial on $\TransInv (T)$.
  Similar statement holds for the invariant $\TransInva$ in 
  $\HFKa$. \qed
\end{cor}

\begin{figure}[ht]
\begin{center}
\setlength{\unitlength}{1mm}
{\includegraphics[height=3cm]{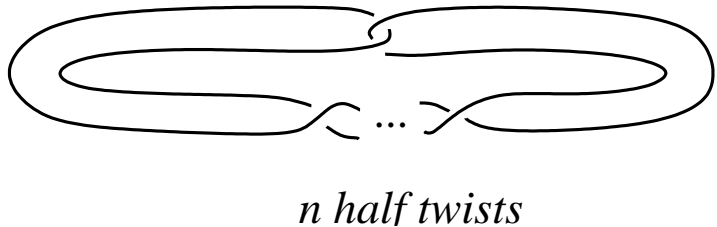}}
\end{center}
\caption{{\bf The Eliashberg--Chekanov knot $E_n$.}}
\label{f:chek}
\end{figure}

This theorem simplifies the computation of $\LegInv$ and $\TransInv $
for many interesting cases, allowing us to use it to distinguish
transversely non--isotopic transverse knots in the same knot type with
the same self--linking number. Recall that a knot type is said to be
{\em transversely simple} if it has no such pairs of transverse
representatives.  The first examples of transversely non--simple knot
types were found by Etnyre-Honda \cite{EHnonsimp} and Birman-Menasco
\cite{BirMen}, and further examples were found in \cite{NOT}.  We will
use Theorem~\ref{t:transz} to show the following:

\begin{thm}\label{t:ecnonsimple}
  The Eliashberg--Chekanov twist knot $E_n$ shown in Figure~\ref{f:chek} is
  not transversely simple for $n$ odd and $n>3$. In fact, for $n$ odd
  there are at least $\lceil\frac{n}{4}\rceil$ transverse knots in the
  standard contact 3--sphere $(S^3 , \xi _{st})$ with self-linking
  number equal to 1, all topologically isotopic to $E_n$, yet not
  pairwise transverse isotopic.
\end{thm}

As a special case, we have the following:

\begin{cor}
The twist knot which is the mirror of $7_2$ in Rolfsen's table is not transversely simple.
\end{cor}
\begin{proof}
The knot $m(7_2)$ is the twist knot $E_5$ from Figure~\ref{f:chek}; hence 
Theorem~\ref{t:ecnonsimple} applies with $n=5$.
\end{proof}

\begin{rem}
  {\rm Recall that by defining and computing the Eliashberg--Chekanov
    DGA for Legendrian representatives of the $5_2$ knot, Chekanov
    \cite{chek} showed that the knot type $5_2$ is not
    \emph{Legendrian simple}, that is, the knot type admits Legendrian
    representatives which are not Legendrian isotopic, though they do
    have the same classical invariants (Thurston--Bennequin and
    rotation numbers). The result was extended in \cite{EFM} for all
    Eliashberg--Chekanov knots $E_n$ $(n\geq 3$ and odd). Notice that
    the DGA's used in these proofs vanish for stabilized knots, hence
    generally these invariants cannot be used to distinguish Legendrian
    approximations of transverse knots. The question of
    whether or not $5_2$ is transversely simple remains open.
    \footnote{In fact, very recently, a classification of the Legendrian and
    transverse isotopy classes of $E_n$ has been announced in~\cite{ENV}}}
\end{rem}

In fact, the same proof leads to finding further transversely non--simple
two-bridge knots in the standard contact 3--sphere; the precise
statement is deferred to Section~\ref{s:ecknots}.  The proof of
Theorem~\ref{t:ecnonsimple}, resting on the transformation rule given
by Theorem~\ref{t:transz}, uses two further ingredients, which we
spell out next. 

Suppose that $F$ is a fixed Seifert surface for $L$. The invariant
$\LegInv (L)\in \HFKm (-Y, L)$ admits an Alexander grading $A_F
(\LegInv (L))$ and (provided the $\SpinC$ structure of the contact
structure has torsion first Chern class) a Maslov grading $M(\LegInv
(L))$.  The Thurston--Bennequin and rotation numbers $\tb (L)$ and
$\rot _F (L)$ of the null--homologous Legendrian knot $L$ can be
defined in the standard way, cf. Section~\ref{s:gradings}.  The
relationship between these numerical invariants of $\LegInv (L)$ and
$L$ is given as follows:

\begin{thm}
  \label{thm:BigradingFormula}
Let $L\subset (Y, \xi )$ be a Legendrian knot in the contact
3--manifold $(Y, \xi )$ and suppose that $F$ is a Seifert surface for
$L$. Then, the chain $\LegInv(L)\in \CFKm(-Y,L)$ is supported in
Alexander grading
  \begin{equation}
    \label{eq:AlexanderGradingFormula}
    2A_F(\LegInv(L))=\tb (L)-\rot _F(L)+1.
  \end{equation}
  If $c_1 (\xi )$ is torsion, then the Maslov grading of $\LegInv (L)$
  is determined by
  \begin{equation}
    \label{eq:MaslovGradingFormula}
    2A(\LegInv(L))-M(\LegInv(L))=d_3(\xi),
  \end{equation}
where $d_3(\xi )$ is the 
3--dimensional invariant (also known as the Hopf invariant) of the
2--plane field underlying the contact structure $\xi $. 
\end{thm}
Since the self--linking number $\sln _F (T)$ of a transverse knot $T$ can
be computed from its Legendrian approximation $L$ as
\[
\sln _F (T)=\tb (L)-\rot _F(L),
\]
we get the following:

\begin{cor}
For a contact 3--manifold $(Y, \xi )$ and transverse knot $T$ the
transverse invariant $\TransInv (T)$ has Alexander grading $A_F
(\TransInv (T))=\frac{1}{2}(\sln _F(T)+1)$ and (provided $c_1(\xi )$
is torsion) Maslov grading $M(\TransInv (T))=\sln(T) +1-d_3(\xi
)$. Similar identities hold for the invariant in
$\HFKa$. \qed
\end{cor}
The second ingredient in the proof of Theorem~\ref{t:ecnonsimple} is a
refinement of the invariant defined in \cite{LOSS}. Recall, that
$\LegInv (L)$ was only defined up to graded automorphisms of the
ambient knot Floer homology group $\HFKm (-Y, L)$. By defining the
action of the mapping class group $\MCG(Y, L)$ of the knot complement
on the knot Floer homology, we will show that the Legendrian isotopy
class of $L$ gives rise to an element in $\HFKm(-Y,L)/\pm \MCG(Y,L)$
rather than in its quotient $\HFKm (-Y, L)/\Aut(Y,L)$ (here, by $\pm
\MCG(Y,L)$ we are emphasizing that we divide out also by the
automorphism gotten by multiplication by $-1$). Since $\pm \MCG(Y,L)$
is typically much smaller than $\Aut(Y,L)$, this lift enables us to
use the invariant much more effectively.  By using $\Field = \Z
/2\Z$--coefficients, the action of $-1$ can be ignored, and in the
following we will apply this choice of coefficient groups.  The
precise formulation of the Heegaard Floer theoretic result showing the
existence of the $\MCG$--action will be given in
Section~\ref{sec:MCG}.

The paper is organized as follows. Section~\ref{s:prelims} is devoted
to the collection of preliminary results, and an explanation how the
invariant is lifted from $\HFKm/\Aut$ to $\HFKm /\pm \MCG$.  The proof
of Theorem~\ref{t:transz} is given in Section~\ref{s:surg}, and we
verify the formulae computing the Alexander and Maslov gradings of
$\LegInv$ in Section~\ref{s:gradings}.  We study Eliashberg--Chekanov
knots --- and certain further two--bridge knots --- in
Section~\ref{s:ecknots}, and in particular give the proof of
Theorem~\ref{t:ecnonsimple}. Finally in Section~\ref{sec:MCG} the
necessary Heegaard Floer theoretic discussion for defining the action
of $\MCG(Y,L)$ on the knot Floer groups is given.

{\bf {Acknowledgements.}} The authors wish to thank Paolo Ghiggini,
Andr{\'a}s Juh{\'a}sz, Paolo Lisca, Zolt\'an Szab\'o and Dylan
Thurston for interesting and helpful discussions. The authors are
especially grateful to Sucharit Sarkar for some discussions of
naturality.  PSO was supported by NSF grant number DMS-0505811 and
FRG-0244663. AS acknowledges support from the Clay Mathematics
Institute.  AS was also partially supported by OTKA T67928 and by
Marie Curie TOK project BudAlgGeo.

\section{Preliminaries}
\label{s:prelims}
We review some of the constructions which will be used throughout the
paper.

\subsection{Knot Floer homology}

We set up notation for knot Floer homology, following~\cite{OSzknot}
(see also~\cite{RasmussenThesis}).  Let $\Sigma$ be a closed, oriented
surface of genus $g$, and $\alphas=\{\alpha_1,...,\alpha_g\}$ be a
$g$--tuple of homologically linearly independent, pairwise disjoint
circles; and let $\betas=\{\beta_1,...,\beta_g\}$ be another such
$g$--tuple of circles.  The triple $(\Sigma , \alphas,\betas)$ is a
Heegaard diagram specifying a closed, oriented $3$--manifold $Y$, built
as follows. We start with the zero-handle, and then regard the
$\alpha$--curves as belt circles of $1$--handles attached to this
zero-handle, and the $\beta$--curves as attaching circles of
$2$--handles. To complete $Y$, we attach the unique $3$-handle.

Fixing two points $z,w\in \Sigma$ in the complement of the $\alpha$--
and $\beta$--curves, an oriented knot $K\subset Y$ is specified as
follows. Connect $z$ to $w$ by a standardly embedded arc disjoint from
the attaching disks in the handlebody determined by the
$\alpha$--curves and $w$ to $z$ by such an arc in the handlebody of
the $\beta$--curves. Notice that the definition, in fact, equips $K$
with an orientation. Consider $\Sym^g(\Sigma)$, equipped with the
totally real submanifolds
\begin{eqnarray*}
\Ta=\alpha_1\times...\times \alpha_g &{\text{and}}&
\Tb=\beta_1\times...\times \beta_g.
\end{eqnarray*}  A suitable adaptation of
Lagrangian Floer homology in this context results in the {\em knot
  Floer homology} groups $\HFKm(Y,K)$, which are the homology groups
of a chain complex $\CFKm(Y,K)$ defined over $\Field[U]$, which is
freely generated by intersection points $\Ta\cap\Tb$.  In the case
where $Y$ is a rational homology three-sphere, these groups are
bigraded,
\[
\HFKm(Y,K)=\oplus _{d\in\Q,\s\in\RelSpinC(Y,K) } \HFKm_d(Y,K,\s),
\]
where $d$ is the {\em Maslov grading} and $\s$, which runs through
{\em relative $\SpinC$ structures} on $Y-K$, is the {\em Alexander
  grading}.  In cases where $Y$ is not a rational homology sphere, we
impose the assumption that $K$ is null-homologous, and we work with
Heegaard diagrams satisfying suitable admissibility hypotheses as
in~\cite{HolDisk}. Even under these hypotheses, when $b_1(Y)>0$, the
Maslov grading is no longer a rational number (except when we consider
relative $\SpinC$ structures whose first Chern class is
torsion). Unless otherwise stated, we will work with $\Field
=\Zmod{2}$ coefficients.

\subsection{Legendrian invariants}
\label{ss:leginv}

We briefly recall the construction of the Legendrian invariant
from~\cite{LOSS} (compare also~\cite{HKM1, OSzcont}).

\begin{figure}[ht]
\begin{center}
\setlength{\unitlength}{1mm}
{\includegraphics[height=4cm]{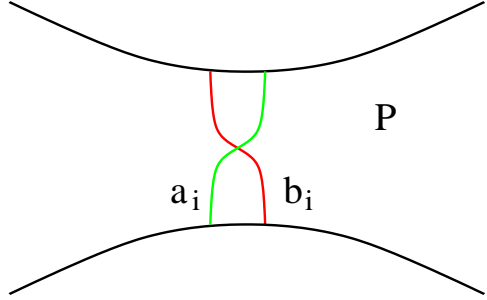}}
\end{center}
\caption{{\bf The arcs $a_i $ and $b_i$.}}
\label{f:ab}
\end{figure}

For a Legendrian knot $L$ in a contact 3--manifold $(Y, \xi )$ a
Heegaard decomposition adapted to this situation can be found by the
following recipe. Choose an open book decomposition of $Y$ compatible
with the contact structure $\xi$ in such a way that $L$ is a
homologically essential curve on one of the pages of the open book.
Such choice is possible, as can be easily verified either by the
application of Giroux's algorithm for constructing open book
decompositions for contact 3--manifolds through contact
cell--decompositions, or by the algorithm of Akbulut and Ozbagci
\cite{AO}, cf. also \cite{Firat}. The open book decomposition, in
turn, provides a Heegaard decomposition, with Heegaard surface given
as the union of two pages $P_{+1}$ and $P_{-1}$ of the open book, and
$\alpha$-- and $\beta$--curves given by the following procedure.
Choose arcs $a_i$ for $i=1,\ldots,n$ in the page $\Page_{+1}$ of the
open book which are disjoint and represent a basis of $H_1(\Page_{+1},
\partial
\Page_{+1})$; i.e., by cutting $\Page_{+1}$ open along the $a_i$ we
get a disk. Let $b_i$ be a slight perturbation of $a_i$, chosen so
that $b_i$ is disjoint from all $a_j$ with $j\neq i$, and intersects
$a_i$ transversely in a single intersection point with orientation
$+1$, as pictured in Figure~\ref{f:ab} (cf. also \cite{HKM1, LOSS}).
In the presence of $L$, the system $\{ a_i \}$ can be chosen in a way
that only $a_1$ intersects $L$, and this intersection is a single
transverse point. Let us take $\alpha _i$ to be the union of $a_i$
with its image under the identity map on the opposite page
$\Page_{-1}$, while $\beta _i$ is the union of $b_i$ together with its
image under the monodromy map $\phi$ (regarded as a map $\phi \colon
P_{+1}\to P_{-1}$) of the open book decomposition.  Clearly,
$\Page_{+1}-a_1-...-a_n-b_1-...-b_n$ consists of $2n+1$
components. For each $i=1,...,n$, we have two components whose
boundary consists of an arc in $a_i$, an arc in $b_i$, and an arc in
$\partial \Page_{+1}$. There is one remaining component (whose
boundary meets all the $a_i$ and $b_i$), and we place the basepoint
$w$ in this region. Moreover, we place $z$ in one of the two remaining
components meeting $a_1$ and $b_1$. We choose this component so that
the induced orientation coincides with the given orientation of the
Legendrian knot $L$. (Recall that we obtain an orientation on $L$ by
orienting its subarc in $\Page_{+1}-a_1-...-a_n$ so as to go from $w$
to $z$.)  In this manner, we obtain a doubly--pointed Heegaard diagram
$(\Sigma,\alphas,\betas,w,z)$ for the oriented Legendrian knot
$L\subset Y$.

\begin{defn}
  \label{def:Adapted}
  The doubly--pointed Heegaard diagram $(\Sigma,\alphas,\betas,w,z)$
  constructed above is said to be {\em adapted} to the open book
  decomposition represented by the surface $\Page$, monodromy map
  $\phi$, and Legendrian knot $L\subset (Y,\xi)$.
\end{defn}

\begin{rem}
{\rm Notice that $z$ and $w$ specify a smooth isotopy class of knots
by connecting $z$ to $w$ on $P$ through $a_1$ and $w$ to $z$ through $b_1$.
In fact, these data uniquely specify a Legendrian isotopy class of 
Legendrian knots, as it is shown in \cite[Theorem~2.7]{LOSS}.}
\end{rem}

Following the convention of \cite{HKM1}, we reverse the roles of the
$\alpha$-- and the $\beta$--circle, hence we examine the Heegaard
decomposition $(\Sigma , \beta , \alpha )$.  Such a change reverses
the orientation of the knot, hence in order to keep the given
orientation of the knot, we also switch the roles of $z$ and $w$,
giving $(\Sigma , \beta , \alpha , z, w)$.  With the reversal of the roles
of the $\alpha$-- and $\beta$--circles, the distinguished
intersection point $\x = (a_1\cap b_1, \ldots , a_n \cap b_n )\in
\Ta\cap \Tb$ becomes  be a cycle in $\CFKm (-Y,L)$. In \cite{LOSS} it was
shown that the homology class $\LegInv (L) \in \HFKm(-Y, L)$
represented by $\x $ is independent of the choices made in its
definition, i.e., from the choice of the adapted open book
decomposition and the basis $\{ a_1, \ldots , a_n \}$.  In addition,
if $L_1, L_2$ are Legendrian isotopic Legendrian knots, then there is
an isomorphism $\HFKm(-Y, L_1)\to \HFKm (-Y, L_2)$ mapping $\LegInv
(L_1)$ to $\LegInv (L_2)$. Viewing $\HFKm (-Y, L)$ as an abstract
group, the element $\LegInv (L_1)$ is therefore defined only up to
identification of $\HFKm (-Y, L_1)$ with $\HFKm (-Y, L_2)$, providing
the following:

\begin{thm}[\cite{LOSS}]\label{t:loss}
  The element $\LegInv(L)$ is an invariant of the oriented Legendrian
  knot $L$, with values in the graded module $\HFKm(-Y,L)$, modulo its
  graded automorphisms. For a transverse knot $T$ and Legendrian
  approximation $L$, the class $\LegInv(L)$ (being invariant under
  negative stabilization) is an invariant 
\[
\TransInv (T)\in \HFKm(-Y,L)/\Aut(\HFKm(-Y,L))
\]
of the transverse isotopy class of $T$.  \qed
\end{thm}

The above theorem supplies an invariant in knot Floer homology, modulo
automorphisms. There are various strengthenings of the above
statements, resulting from the types of restrictions one can naturally
place on the allowed automorphisms. An example of such strengthening
was given in \cite{LOSS} for connected sums. In a slightly different
direction, on the crudest level, if $T_1$ and $T_2$ are two different
transverse realizations of the same knot, and $\TransInv$ for one of
them vanishes while for the other does not, then the above theorem
ensures that $T_1$ and $T_2$ are not transversely isotopic. But
$\HFKm$ has more algebraic structure than merely a bigraded
$\Field[U]$--module: it is naturally the homology group of an associated
graded object of a filtered complex; as such, it comes equipped with
higher differentials. Thus, if $T_1$ and $T_2$ both have {\em
non--vanishing} transverse invariant, but $d_1(\TransInv(T_1))$
vanishes while $d_1(\TransInv(T_2))$ does not, then $T_1$ and $T_2$
are not transversely isotopic. This refined structure was used to
distinguish transversally non--isotopic knots in the combinatorial
context in~\cite{NOT}. 

Sometimes such algebraic properties are insufficient to distinguish
transverse knots, and it becomes necessary to use more refined
geometric tools. Below we will describe the lift of the Legendrian
invariant from $\HFKm /\Aut$ to $\HFKm /\pm \MCG$. To this end,
consider a knot $K\subset Y$ and let $\DiffP(Y,K)$ denote the space of
diffeomorphisms from $Y$ to itself which fix $K$ pointwise (or, equivalently, fix a point $p$ on $K$). Let $\DiffP_0(Y,K)$
be the set of those elements in $\DiffP(Y,K)$ which can be connected
to the identity map (through a one--parameter family of elements in
$\DiffP(Y,K)$). Let $\MCG(Y,K)$ denote the \emph{mapping class group}
of a knot complement, that is,
\[
\MCG(Y,K)=\frac{\DiffP(Y,K)}{\DiffP_0(Y,K)}.
\]
The tools of~\cite{HolDiskFour}, adapted to the context of links, lead
to an induced action of $\MCG(Y,K)$ on $\HFKm(Y,K)$, which will be
spelled out in Section~\ref{sec:MCG}. More generally, a diffeomorphism
of $(Y,K,p)$ to $(Y',K',p')$ (where here $p\in K\subset Y$ and $p'\in
K'\subset Y'$) induces a well-defined map on $\HFKm$. (See
Theorem~\ref{thm:Naturality}, and the remarks afterwards.)  This
concept leads us to the following refinement of Theorem~\ref{t:loss}:

\begin{thm}
  If $L\subset (Y,\xi)$ is a null--homologous Legendrian knot
  resp. $T\subset (Y, \xi )$ is a null--homologous transverse knot,
  then the invariant $\LegInv(L)$, resp. $\TransInv(T)$ naturally
  takes values in $\HFKm(-Y,L)/\pm \MCG(Y,L)$; i.e. if $L_1$ and $L_2$
  are Legendrian resp. transverse realizations of the same knot type
  $K$ whose invariants $\LegInv$ resp. $\TransInv$ lie in different
  orbits in $\HFKm(-Y,K)$ under the group generated by multiplication
  by $-1$ and the mapping class group action of the knot complement,
  then $L_1$ and $L_2$ are not Legendrian resp. transversely isotopic.
\end{thm}
\begin{proof}
Fix a knot $L\subset Y$ in the knot type $K$, and fix a point $p\in L$,
 and consider $\HFKm(-Y,L,p)$.  For a Legendrian representative $L_1$ of
$K$ and a point $p_1\in L_1$,
consider an isotopy $\varphi _t$ between $(L,p)$ and $(L_1,p_1)$ with
time-one map $\varphi _1$ inducing the isomorphism $(\varphi _1)_*
\colon \HFKm (-Y,L_1,p_1)\to \HFKm (-Y, L,p)$ on the knot Floer
homologies. Consider the image of $\LegInv (L_1)$ (defined up to a
multiplication by $(-1)$) in $\HFKm(-Y, L,p)$; this element will depend
on the chosen isotopy $\varphi _t$. Another isotopy $\psi _t$ will
give rise to another identification $(\psi _1)_*$, for which the
composition $(\psi _1)_* \circ (\varphi _1)_*^{-1}$ is the action of
the mapping class $\psi _1\circ \varphi _1 ^{-1}\in \MCG(-Y,L)$ on
$\HFKm (-Y, L,p)$. Therefore the element $(\varphi _1)_*(\LegInv (L_1))$
is well--defined up to the action of $\pm \MCG(-Y, L)$.

Suppose that $(L_2,p_2)$ is another Legendrian knot with the property
that $L_1$ and $L_2$ are Legendrian isotopic, and let $\zeta_t$ be the
Legendrian isotopy between $L_1$ and $L_2$. Note that for any
Legendrian knot $L_0$ and $p,p'\in L_0$, $(L_0,p)$ and $(L_0,p')$ are
Legendrian isotopic (as can be verified using contact Hamiltonian
functions). Thus, we can assume that the time-one map $\zeta_1$ carries 
$p_1$ to $p_2$. Therefore, we have an induced map
$(\zeta_1)_* \colon \HFKm (-Y, L_1)\to \HFKm (-Y, L_2)$, which by
\cite[Corollary~3.6]{LOSS} maps $\LegInv (L_1)$ to $\LegInv
(L_2)$. Hence the composition of the isotopies shows that the $\pm
\MCG(-Y,L)$--orbit of the image of $\LegInv (L_1)$ is equal to the
$\pm \MCG(-Y,L)$--orbit of the image of $\LegInv (L_2)$, concluding
the proof. When using $\Z / 2\Z$--coefficients, multiplication by
$(-1)$ induces the trivial action, hence $\LegInv (L)$ is defined as
an $\MCG (Y,L)$--orbit in this case.
\end{proof}

In the following (in order to keep the discussion simpler) we will
still refer to $\LegInv $ as an element of the knot Floer homology
$\HFKm$, although the particular element is just a representative of
the corresponding $\pm \MCG$--orbit. Some of the basic properties of
$\LegInv$ from~\cite{LOSS} are summarized in the following:

\begin{thm}[\cite{LOSS}]
If the contact invariant $c(Y, \xi )$ is nonzero, then any Legendrian knot
$L\subset (Y, \xi ) $ has nonvanishing $\LegInv$--invariants. If 
$c(Y, \xi )=0$ then $\LegInv(L)$ is a $U$--torsion class. If $(Y, \xi )$
is overtwisted and $L$ is a loose knot (that is, $Y-L$ is overtwisted) then
$\LegInv(L)=0$. \qed
\end{thm}

In fact, the nonvanishing result was shown by applying 
the map 
\[
\HFKm (-Y, K, \s )\to \HFa (-Y, \s)
\]
given by the specialization $U=1$, which maps the Legendrian invariant
$\LegInv(L)$ to the contact invariant $c(Y,\xi )$ of the contact
structure $(Y, \xi )$ of~\cite{OSzcont}.  In \cite{LOSS} a number of
explicit computations for $\LegInv(L)$ were given by choosing the
appropriate open book decompositions, and determining the homology
class of the intersection point $\x$ from a direct analysis of the
chain complex.  In the next section we will show another way of
computing $\LegInv$, which now will rely on a transformation rule
developed for contact $(+1)$--surgeries. In this argument we will need
to understand how knot Floer homology behaves under a map associated
to a surgery.

\subsection{Maps induced by surgery}
\label{ss:mapsbysurgery}

Suppose that $Y$ is a three-manifold equipped with a framed knot
$C$. Let $Y_f(C)$ denote the three-manifold obtained as 
surgery with the prescribed framing $f$ along $C$ in $Y$.  The triple
$(Y, C, f)$ can be described by a \emph{Heegaard triple} $(\Sigma
,\alpha , \beta , \gamma )$, where $(\Sigma ,\alpha , \beta )$ gives
$Y$ and $(\Sigma ,\alpha , \gamma )$ gives $Y_f (C)$. By counting
holomorphic triangles in $\Sym ^g (\Sigma )$ with boundaries on the
totally real tori ${\mathbb {T}} _{\alpha}, {\mathbb {T}}_{\beta }$
and ${\mathbb {T}} _{\gamma }$, and choosing a particular cycle in the
chain complex of $(\Sigma , \beta , \gamma )$, we get a map
\[
{\widehat {F}}_C\colon \HFa (Y) \to \HFa (Y _f (C)).
\]
(It is shown in \cite{HolDiskFour} that the map ${\widehat {F}}_C$
does not depend on the particular choices and, in fact, is an
invariant of the 4--dimensional surgery cobordism.)

If $\mu$ denotes a meridian for $C$, then we can think of $f+\mu$ as a new framing. In fact, $Y_{f+\mu
}(C)$ can be regarded as surgery along a framed knot $(C',f')\subset Y_f
(C)$ and $Y$ can be regarded as the result of a surgery along $(C'',
f'')\subset Y_{f+\mu }(C)$.  If $Y_1$, $Y_2$, and $Y_3$ are three
three--manifolds which are related in this manner, i.e. $Y_1=Y$, $Y_2=
Y_f(C)$, and $Y_3=Y_{f+\mu }(C)$, then we say that the cyclically
ordered triple $(Y_1,Y_2,Y_3)$ forms a {\em distinguished triangle}.
Note that the roles of $Y$, $Y_f(C)$, and $Y_{f+\mu }(C)$ are
cyclically symmetric. In fact, all three three-manifolds are obtained
by Dehn filling the same three-manifold $M$ with torus boundary along
three different surgery slopes, which meet pairwise in a single point.
According to~\cite{HolDiskTwo}, the maps ${\widehat {F}}_C$,
${\widehat {F}}_{C'}$ and ${\widehat {F}}_{C''}$ fit into an exact
triangle.

This construction can be refined to the case of knot Floer homology,
as in~\cite{OSzknot}. Specifically, suppose that $K\subset Y-C$ is a
null--homologous knot. In this
case, $K$ naturally induces a null--homologous knot $K'$ inside
$Y_f(C)$, and $C$ gives rise to the map 
\[
{\widehat F}_{C}\colon \HFKa(Y,K) \longrightarrow \HFKa(Y_f(C),K') 
\]
on knot Floer homology groups. Since the assumption of $K\subset Y-C$
being null--homologous implies that the linking number of $K$ with $C$
is trivial, the induced map also preserves Alexander grading,
see~\cite[Proposition~8.1]{OSzknot}.

\begin{defn}
  Suppose that $C\subset Y$ is a framed knot and $K \subset Y-C$ is a
  null--homologous knot (i.e. one whose linking number with $C$ is
  trivial).  $K$ can be thought of as a knot in $Y$, in $Y_f(C)$
  (denoted by $K'$), or in $Y_{f+\mu }(C)$ (denoted by $K''$).  We
  call the triple of knots $\{(Y,K), (Y_f(C),K'), (Y_{f+\mu
  }(C),K'')\}$ a {\em distinguished triangle of knots}.
\end{defn}

Note once again that the roles of the three knots in a distinguished
triangle are cyclically symmetric. The surgery exact triangle 
of \cite{HolDiskTwo} then has the following extension:

\begin{thm}
  \cite[Theorem~8.2]{OSzknot} If
    $\{(Y_1,K_1),(Y_2,K_2),(Y_3,K_3)\}$ is a distinguished triange of
    knots, then the corresponding Alexander grading
    preserving maps fit into an exact triangle
\vskip1in
\begin{center}
\begin{graph}(5,-4)
\graphlinecolour{1}\grapharrowtype{2}
\textnode {A}(1,1.5){$\HFKa (Y_1,K_1)$}
\textnode {B}(5, 1.5){$\HFKa (Y_2,K_2)$}
\textnode {C}(3, 0){$\HFKa (Y_3,K_3)$}
\diredge {A}{B}[\graphlinecolour{0}]
\diredge {B}{C}[\graphlinecolour{0}]
\diredge {C}{A}[\graphlinecolour{0}]
\freetext (3,1.8){$\Fa_ 1$}
\freetext (4.6,0.6){$\Fa _2$}
\freetext (1.4,0.6){$\Fa _3$}
\end{graph}
\end{center}
\end{thm}

\begin{rem}
{\rm Notice that the independence of the maps ${\widehat {F}}_i$ on
  the chosen Heegaard triple is not claimed above --- although it is
  plausible to expect that these maps will depend only on the
  4--dimensional cobordism defined by the surgery, cf. also
  \cite{Sahamie}.  In our arguments we will use only the construction
  of these maps (given by counting holomorphic triangles) and the
  exactness of the triangle associated to a distinguished triangle.}
\end{rem}

The map ${\widehat F}_C$ on $\HFKa$ defined by counting
pseudo-holomorphic triangles can be extended to the case of $\HFKm$ as
well, to define a map
\[
F_{C,\spinc}\colon \HFKm (Y ) \to \HFKm (Y_f (C)).
\]
Just as in the case for closed three-manifolds, in order for this map
to be well-defined, we must fix a $\SpinC$-equivalence class of
pseudo-holomorphic triangle on the surgery cobordism.  For the case of
$\HFKa$, this choice can be omitted, with the understanding that
${\widehat F}_C$ is obtained as a sum over all such choices. In the
case of $\HFKa$ this is allowed since the map on $\CFa$ is trivial for
all but finitely many $\SpinC$-equivalence classes of maps
(compare~\cite[Theorem~3.3]{HolDiskFour}); but again just as in the
closed case, this is no longer the case for $\HFKm$, and a choice of a
Spin$^c$ structure $\spinc$ on the 4--dimensional surgery cobordism is
necessary.

\section{The effect of contact $(+1)$--surgery on the Legendrian invariant}
\label{s:surg}

Consider now a null--homologous Legendrian knot $L\subset (Y, \xi)$
and another Legendrian knot $S\subset (Y, \xi )$ (disjoint from $L$),
and perform Legendrian $(+1)$--surgery along $S$. The resulting
contact 3--manifold $(Y_S, \xi _S)$ obviously contains $L$ as a
Legendrian knot, denoted by $L_S$; suppose that $L_S$ is still
null--homologous in $Y_S$. (This condition holds, for example, if the
linking number of $L$ and $S$ is zero.) Let us choose an open book
decomposition which is adapted to $(Y, \xi , S, L)$ as in
Definition~\ref{def:Adapted}; that is, the open book decomposition is
compatible with the contact structure $\xi $ and contains $S,L$ on one
of its pages.  We can further assume that $S$ and $L$ are
homologically essential in the page $\Page$, represent different
homology elements, and the complement of $S$ in $\Page$ is connected.

The open book decomposition gives rise to a Heegaard diagram
$(\Sigma,\alphas,\betas,w,z)$: the Heegaard surface $\Sigma$ is the
union of two pages $\Page_{+1}$ and $\Page_{-1}$, which are oriented
oppositely and the $\alpha$-- and $\beta$--curves on $\Page_{+1}\cup
(-\Page_{-1})$ are defined as $\alpha _i= a_i\cup {\overline a_i}$ and
$\beta _i=b_i \cup {\overline{\phi (b_i)}}$, where $\phi$ represents
the monodromy of the open book, $\{ a_i\}$ is a basis in $\Page_{+1}$
intersecting $L$ in a unique point, and $x\mapsto{\overline x}$
represents the (orientation-reversing) identification between
$\Page_{+1}$ and $(-\Page_{-1})$. Recall that the Legendrian invariant
$\LegInv(L)$ is represented by the intersection point
$$\x_0=\prod_{i=1}^g (b_i\cap a_i) \in\Tb\cap\Ta,$$
now thought of as
a generator for $\CFKm(-Y,L)$, specified by the Heegaard diagram
$(\Sigma,\betas,\alphas,z,w)$.

Since $L$ and $S$ are distinct in homology, we can assume that $a_1$
intersects $L$ transversely in a unique point (and it is disjoint from
$S$), $a_2$ intersects $S$ transversely in a unique point (and it is
disjoint from $L$), and $a_i$ for $i>2$ are disjoint from both $L$ and
$S$.  By our choice, the open book decomposition gives a compatible
open book for $(Y_S, \xi _S, L_S )$ as well: just compose the
monodromy $\phi$ with the left--handed Dehn twist $D^{-1}_S$
(corresponding to the fact that we perform contact $(+1)$--surgery)
along $S$. This results in a change of the Heegaard diagram for
$(Y,L)$ by applying the Dehn twist $D_S$ to all $\phi (b_i)$
intersecting $S$ on $P_{-1}$. In the resulting diagram it is rather
complicated to detect the effect of the handle attachment.

To simplify matters, we set up things slightly differently, as
follows.  Consider the doubly--pointed Heegaard diagram
$(\Sigma,\alphas,\gammas,w,z)$, where here the basepoints $w$ and $z$
are placed adjacent to $a_1$ and $b_1$ as before.  The curve
$\gamma_2$ is a small perturbation of
$D_S(b_2)\cup{\overline{\phi(b_2)}}$, while for $i\neq 2$, $\gamma_i$
is a small perturbation of $\beta_i$.
Although $(\Sigma,\alphas,\gammas,w,z)$ does represent $(Y_S,L_S)$,
it is  {\em not} adapted to 
the Legendrian knot $L_S \subset (Y_S, \xi _S)$, in the sense of
Defintion~\ref{def:Adapted}. However, if we consider the Heegaard
diagram $(\Sigma,\deltas,\gammas,w,z)$ where here $\delta_i$ are
suitable small perturbations of the $\alpha_i$ for all $i\neq 2$, while
$\delta_2$ is a perturbation of $D_S(a_2)\cup {\overline {D_S(a_2)}}$,
then it is easy to see that
$(\Sigma,\deltas,\gammas,w,z)$ is an adapted Heegaard diagram for
$(Y_S,\xi_S, L_S)$, in the sense of Definition~\ref{def:Adapted}.

The invariant $\LegInv(L)\in\CFKm(-Y,L)$ is represented by a cycle
$\x_0\in\Tb\cap\Ta$ for the Heegaard diagram
$(\Sigma,\betas,\alphas,z,w)$; similarly, $\LegInv(L_S) $ is
represented by the intersection point point $\x_1\in\Tc\cap\Td$ for
the Heegaard diagram $(\Sigma,\gammas,\deltas,z,w)$, thought of as an
element of $\CFKm(-Y_S,L_S)$.  In order to relate the Legendrian
invariants $\LegInv(L)$ and $\LegInv(L_S)$ we will use the
intermediate diagram $(\Sigma, \gammas,\alphas,z,w)$
for $(-Y_S,L_S)$.

\begin{figure}[ht]
\mbox{\vbox{\epsfbox{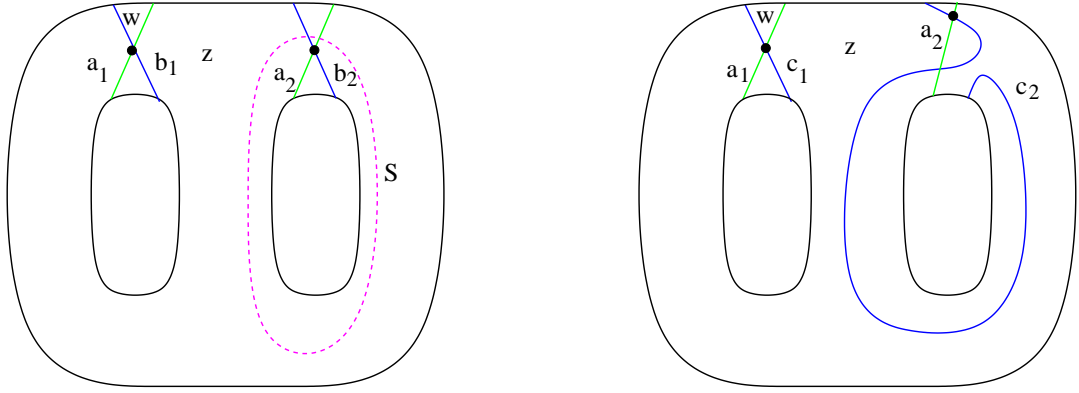}}}
\caption{\label{f:IntermediateCycle} {\bf{Intermediate cycle.}}  We
  start from an open book as on the left, which supports both our
  initial Legendrian knot $L$, and also the Legendrian knot $S$ along
  which surgery is to be performed (note that this diagram is, in
  general, stabilized by the addition of further one--handles equipped
  with $a_i$-- and $b_i$--curves).  $S$ here is represented by the
  dotted line, and $L$ is recorded by the pair of basepoints $w$ and
  $z$. The curve $c_2=D_{S}(b_2)$ on the right is obtained by Dehn
  twist along $S$ (after introducing a finger move).  The components
  of the ``intermediate intersection point'' $\y$ are indicated by the
  dots on the right--hand side diagram. (Recall that the Heegaard
  diagram $(\Sigma,\gammas,\alphas,z,w)$ for $(-Y_S,L_S)$ is gotten by
  $\gamma_i=c_i\cup {\overline{\phi(b_i)}}$, $\alpha_i=a_i\cup
  {\overline {a_i}}$.)}
\end{figure}

Specifically, we would like to find an intersection point
$\y\in\Tc\cap\Ta$ representing the Legendrian invariant for knot Floer
homology of $(-Y_S,L_S)$, only using the Heegaard diagram
$(\Sigma,\gammas,\alphas,z,w)$. In this case $\gamma_i$ and $\alpha_i$
meet in a single point $x_i$ in $P_{+1}$ for all $i\neq 2$.  Special
care must be taken for $i=2$. Recall that $\gamma_2$ and
$\alpha_2$ on $\Page_{+1}$ are a perturbation $c_2$ of $D_{S}(b_2)$,
and $a_2$ respectively. If care is not taken, these curves will be
disjoint on $\Page_{+1}$.  However, we make a finger move on $c_2$ to
ensure it meets $a_2$ in two points, as pictured in
Figure~\ref{f:IntermediateCycle}, creating an intersection point on
$P_{+1}$ representing the Legendrian invariant. Now,
$\gamma_2\cap\alpha_2$ consists of two points $y_1$ and $y_2$ on
$\Page_{+1}$.  For one of these choices $y_2$, we have that $\y
=(x_1,y_2, x_3, \ldots ,x_n)$ is a cycle in $\Tc\cap\Ta$ (where
the other components $x_i$ are $a_i\cap b_i$ as before), since there
is no positive domain $D$ supported in $\Page_{+1}$ with
$n_{z}(D)=0$ and with initial point $\y$. This choice is illustrated
in Figure~\ref{f:IntermediateCycle}.

\begin{lem}
  \label{l:IntermediateCycle}
  The intersection point $\y\in\Tc\cap\Tb$ $($as represented by the
  diagram $(\Sigma,\gammas,\alphas,z,w))$ is a cycle in
  $\HFKm(-Y_S,L_S)$.
\end{lem}

\begin{proof}
  The statement is proved using the same mechanism which guarantees
  that the distinguished intersection point for an adapted Heegaard
  diagram is a cycle, cf.~\cite{LOSS,HKM1}. Specifically, we argue
  that there is no nontrivial class $\phi\in\pi_2(\y,\y')$ with
  $n_z(\phi)=0$.  This follows from a direct analysis of the Heegaard
  diagram.
\end{proof}

Recall now that
$$F_{S,\spinc}\colon \CFKm(-Y,L) \longrightarrow \CFKm(-Y_S,L_S)$$
is the map induced by the map
$$\CFKm(\Sigma,\gammas,\betas,z,w)\otimes
\CFKm(\Sigma,\betas,\alphas,z,w)\longrightarrow
\CFKm(\Sigma,\gammas,\alphas,z,w),$$
which is obtained by counting holomorphic triangles representing
the $\SpinC$ structure $\spinc$, after
pairing with the canonical top--dimensional homology class $\Theta$ for
$\CFKm(\Sigma,\gammas,\betas,z,w)$. (We can think of $\Theta$ as represented
by some intersection point between $\Tb$ and $\Tc$.)

\begin{lem}\label{l:elso}
  There is a unique $\SpinC$ structure $\spinc$ on the cobordism 
  from $Y$ to $Y_S$ for which the induced map
  $F_{S,\spinc}(\x_0)$ is nontrivial, where here
  $\x_0=\x (Y, \xi , L)$ represents $\LegInv(Y,
  \xi , L)$. For that choice, we have
  that 
  $$F_{S,\spinc}(\x_0)=\y =(x_1,y_2, x_3, \ldots , x_n).$$
\end{lem}

\begin{proof}
  It is easy to see that the top--dimensional homology class of 
  $$\HFm (\#^{g-1}(S^1\times S^2))=\HFKm (\Sigma , \gamma , \beta , z,
  w)$$ is represented by an intersection point $\Tc\cap\Ta$ supported
  on the page $\Page_{+1}$, which we denote by
  $\Theta_{\gamma\alpha}$. Moreover, there is a plainly visible
  Whitney triangle $\psi_0\in\pi_2(\Theta_{\gamma\alpha},\x_0,\y)$, as
  illustrated on the left-hand picture in
  Figure~\ref{f:IntermediateTriangles}. By the Riemann mapping
  theorem, this triangle has a unique pseudo--holomorphic
  representative.  Let $\spinc$ denote the $\SpinC$ structure induced
  by this pseudo--holomorphic triangle.

\begin{figure}[ht]
\mbox{\vbox{\epsfbox{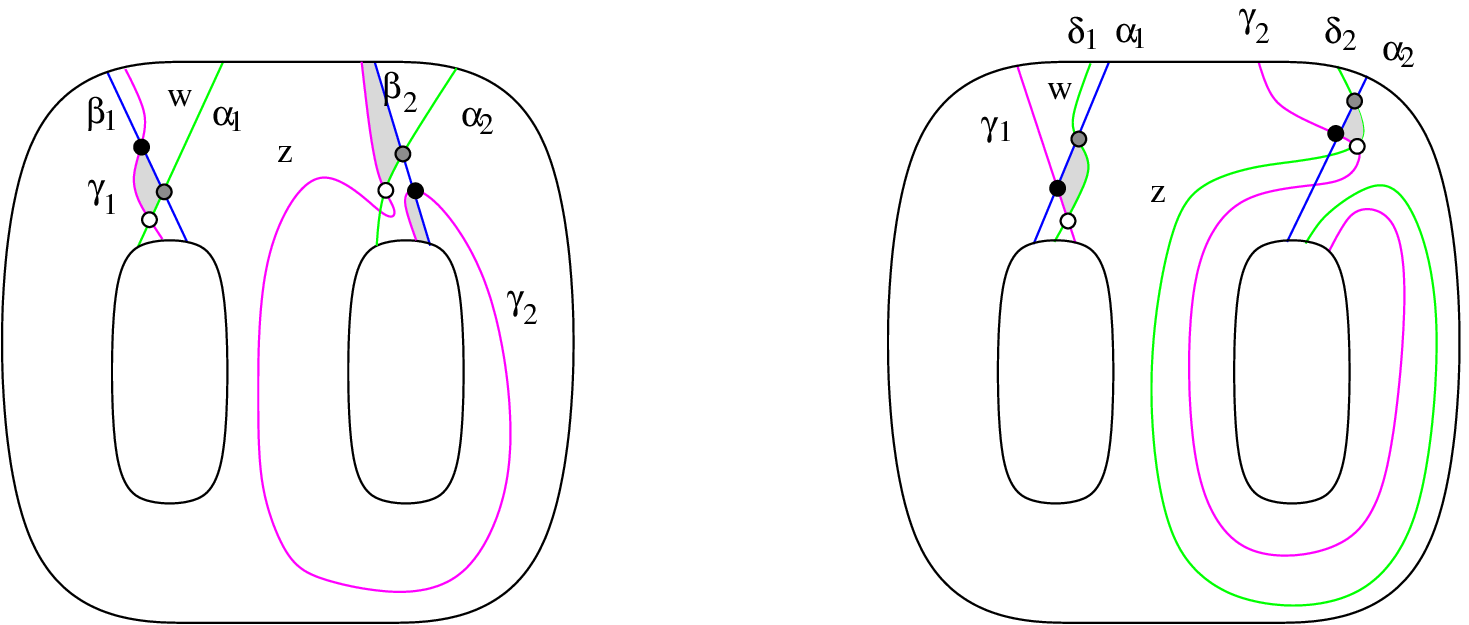}}}
\caption{\label{f:IntermediateTriangles}
  {\bf{Intermediate triangles.}}  On the left, we have illustrated the
  triangle in the Heegaard triple
  $(\Sigma,\gammas,\betas,\alphas,z,w)$, and on the right, we have
  triangles for the Heegaard triple
  $(\Sigma,\gammas,\alphas,\deltas,z,w)$.  The dark circles represent
  initial intersection points ($\Theta_{\gamma\beta}\in\Tc\cap\Tb$ on
  the left, $\y\in\Tc\cap\Ta$ on the right), gray circles represent
  intermediate ones ($\x_0\in\Tb\cap\Ta$ on the left, and
  $\Theta_{\alpha\delta}\in\Ta\cap\Td$ on the right), while the light
  ones represent the final intersection points ($\y\in\Tc\cap\Ta$ on
  the left, $\x_1\in\Tc\cap\Td$ on the right).}
\end{figure}

  We claim that if $\psi\in\pi_2(\Theta_{\gamma\alpha},\x_0,\x_2)$ is
  any homotopy class of Whitney triangles with positive underlying
  domain and $n_{z}(\psi)=0$, then $\x_2=\y$ and $\psi=\psi_0$. To see
  this, we argue that any such $\psi$ has the form $\psi_0*\phi$ for
  some $\phi\in\pi_2(\x_2,\y)$ with $n_{z}(\phi)=0$. As in the proof
  of Lemma~\ref{l:IntermediateCycle}, $\phi$ must be either trivial
  (the case where $\psi=\psi_0$) or it must have a negative local
  multiplicity somewhere. In the latter case, it is easy to see that
  $\psi_0 * \phi$ must also have a negative local multiplicity
  somewhere.  It now follows that $F_{S, \spinct } (\x _0 )=0$ for all
  $\spinct \neq \spinc$, and $F_{S, \spinc }$ maps $\LegInv(Y,\xi,L)$
  to $\y$, as claimed.
\end{proof}

Next we will show that $\y$, in fact, represents the Legendrian
invariant $\LegInv(Y_S, \xi _S, L_S)$.

\begin{lem}\label{l:masodik} 
  Under the homotopy equivalence 
  \[
  H\colon \CFKm (\Sigma,\gammas,\alphas,z,w) \longrightarrow
  \CFKm (\Sigma,\gammas,\deltas,z,w)
  \]
  given by handleslides, the intersection point $\y$ is mapped to
  $\x_1$, representing the Legendrian invariant of $L_S$.
\end{lem}

\begin{proof}
  The appropriate homotopy equivalence is now induced by counting
  pseudo--holomorphic triangles in the Heegaard triple
  $(\Sigma,\beta,\gamma,\delta,z,w)$.  As indicated on the right--hand
  picture in Figure~\ref{f:IntermediateTriangles}, we can find an
  intersection point $\Theta_{\alpha\delta}\in\Ta\cap\Td$ which
  represents the top--most homology class of 
\[
\HFm(\#^{g}(S^2\times S^1))=\HFKm (\Sigma , \alpha, \delta , z,w). 
\]
Once again, we have a Whitney triangle,
  $\psi_0\in\pi_2(\y,\Theta_{\alpha\delta},\x_1)$, evident on the page
  $P_{+1}$, also illustrated on the right--hand picture in
  Figure~\ref{f:IntermediateTriangles}, where $\x_1\in\Tc\cap\Td$
  represents the Legendrian invariant. As before, if
  $\psi\in\pi_2(\y,\Theta_{\alpha\delta},\x')$ is any homotopy class
  with positive domain and $n_z(\psi)=0$, then $\x'=\x_1$ and
  $\psi=\psi_0$.  Thus, $H(\y)=\x_1$ as claimed.
\end{proof}

\begin{proof}[Proof of Theorem~\ref{t:transz}]
  According to Lemma~\ref{l:elso} the map $F_{S,\spinc}(\LegInv(L))$
  is nontrivial for only one choice of $\spinc$, for which it maps
  $\LegInv( L)$ to a homology class represented by $\y$ (which is a
  cycle, according to Lemma~\ref{l:IntermediateCycle}), and by
  Lemma~\ref{l:masodik} the cycle $\y$, in fact, represents
  $\LegInv(Y_S, \xi _S, L_S)$, concluding the proof.
\end{proof}

In applications, it is sometimes more convenient to work with $\HFKa$,
the specialization of $\HFKm$ to $U=0$. Specifically, recall that the
specialization $U=0$ provides a map
\[
\HFKm (-Y, L) \to \HFKa (-Y, L),
\]
and the image of the Legendrian invariant $\LegInv (L)$ under this map
is denoted by $\LegInva (L)$. We have the corresponding maps $\Fa
_{S,\spinc}$ induced by counting pseudo-holomorphic triangles in the
$U=0$ context. In fact,  by writing
$$\Fa_{S}=\sum_{\spinc\in\SpinC(W)} \Fa_{S,\spinc},$$
Theorem~\ref{t:transz} has the following specialization:

\begin{cor}
  Suppose that $(Y_S, \xi _S, L_S )$ is the result of contact
  $(+1)$--surgery along $S\subset (Y, \xi , L)$ and suppose futhermore
  that $L\subset Y$ and $L_S \subset Y_S$ are both null-homologous
  Legendrian knots. Then, for the induced map $F_S$ on $\HFKa(-Y, L)$
  we have
\[
\Fa_{S}(\LegInva (Y, \xi , L))=\LegInva (Y_S, \xi _S, L_S).
\]
\qed
\end{cor}

\section{Gradings}
\label{s:gradings}
\newcommand\OneHalf{\frac{1}{2}}
\newcommand\relspinc{\mathfrak{t}}
\newcommand\cm{\cdot}

In this section we turn our attention to the proof of
Theorem~\ref{thm:BigradingFormula}, relating the bigrading of the
Legendrian invariant of a Legendrian knot $L$ with the classical
Legendrian invariants of $L$.  We start our discussion when the
ambient 3--manifold $Y$ is an integral homology 3--sphere, and
consider the general case at the end of the subsection.

\subsection{Alexander gradings}

  Suppose that $L\subset Y$ is an oriented knot in the intergral
  homology sphere $Y$. Let $F$ denote a Seifert surface for $L$.
  There is a natural map $\RelSpinC(Y,L) \longrightarrow \Z$ given by
\begin{equation}\label{e:Algra} 
\relspinc\mapsto  \OneHalf \langle c_1(\relspinc), F\rangle ,
\end{equation}
 where a \emph{relative $\SpinC$ structure} $\relspinc \in
\RelSpinC(Y,L)$ is regarded as a $\SpinC$ structure ${\hat {\relspinc
}}$ on the $0$--surgery along $L$. The pairing $\langle c_1(\relspinc),
F\rangle$ is interpreted as integration in the result of the
$0$--surgery, i.e. 
\[
\langle c_1(\relspinc), F\rangle =\langle c_1({\hat {\relspinc}}),
{\hat {F}}\rangle ,
\]
where ${\hat {F}}$ is the surface we get by capping off the Seifert
surface $F$. Since $Y$ is an integral homology sphere, the result will
be independent of the particular choice of $F$.

Since any intersection point $\x \in {\mathbb {T}} _{\alpha } \cap
{\mathbb {T}} _{\beta}$ for $(\Sigma , \beta , \alpha , z, w )$
determines a relative $\SpinC$ structure, in view of the above
definition we have an integral--valued Alexander grading belonging to
each intersection point in $\Ta\cap\Tb$.

\begin{thm}
  \label{thm:SpecialCaseAlex}
  The integral Alexander grading of the Legendrian invariant $\LegInv
  (L)$ is given by the formula
  $$2A(\LegInv(L))=\tb(L)-\rot(L)+1.$$
\end{thm}

We start with some basic algebraic topology for open book
decompositions, and the corresponding interpretation of the rotation
number. Then, we turn to an interpretation of the rotation number and
the Thurston--Bennequin invariant in terms of a compatible Heegaard
diagram.  This will lead to a proof of
Theorem~\ref{thm:SpecialCaseAlex}.

Recall that the open book decomposition can be given as a
surface-with-boundary $P$, together with a mapping class $\phi\colon P
\longrightarrow P$ (fixing the boundary). We can form the mapping
torus, which is a three--manifold with torus boundary. Filling the
tori (with appropriate slope), we form a three--manifold $Y(\phi)$
equipped with an open book decomposition. By applying positive
stabilizations, we can assume that the binding $\partial P$ is
connected; we will always assume this extra hypothesis on our chosen
open book decomposition.

\begin{lem}
  \label{lemma:NullHomologousKnot}
  An element $[L]\in H_1(P)$ is in the kernel of $H_1(P)\longrightarrow
  H_1(Y(\phi))$ if and only if it can be written as $[L]=\phi_*(Z)-Z$
  for some $Z\in H_1(P)$. 
\end{lem}

\begin{proof}
  Note that $L$ has linking number zero with the binding. It then
  follows that if the homology class $[L]$ is in the kernel of
  $H_1(P)\longrightarrow H_1(Y(\phi))$, then $L$ is already
  null--homologous in the mapping torus $M(\phi)$ of $\phi$.  Now,
  $M(\phi)$ is homotopy equivalent to the two-complex obtained from
  $P$ by attaching cylinders whose boundary components have the form
  $\phi _*(Z)-Z$. The result now follows.
\end{proof}

\begin{defn}
  Let $f\colon \Pd \longrightarrow \Sigma$ be a map of a two-manifold
  $\Pd$ with boundary into a surface $\Sigma$, which has the property
  that the boundary of $\Pd$ is immersed in $\Sigma$. The {\em Euler
  measure} of this map is defined as the relative Chern number of
  $f^*(\Sigma)$, relative to natural trivialization of its boundary
  inherited from the boundary, thought of as immersed curves in
  $\Sigma$. The Euler measure depends on $f$ only through its induced
  underlying two--chain. For more on the Euler measure,
  see~\cite[Section~7.1]{HolDiskTwo}. The Euler measure of $\Pd$
  will be denoted by $e(\Pd)$.
\end{defn}

We descibe a way to identify the rotation number of the Legendrian
knot in the Heegaard diagram.  The contact distribution defines a
complex line bundle over $Y(\phi)$ whose restriction to the the
mapping torus (thought of as a subset of $Y(\phi)$) coincides with its
``tangents along the fiber''.

\begin{lem}
  \label{lemma:RotationNumber}
  Let $L\subset Y(\phi)$ be a null-homologous knot supported in a
  fiber $P$ for the open book decomposition. Let $\Pd$ be a two-chain
  with $\partial \Pd = [L]+ (Z-\phi_*(Z))$ for the one--cycle $Z$
  found in Lemma~\ref{lemma:NullHomologousKnot}. Then, the rotation
  number of the Legendrian knot $L$ is calculated by the Euler measure
  of $\Pd$.
\end{lem}

\begin{proof}
  We can construct a two--chain $F$ with boundary $L$ as follows. $F$
  is composed of $\Pd$, thought of as supported in a fiber of the open
  book, and then along each boundary cycle of type $Z$, we attach a
  cylinder which traverses the mapping torus, meeting $P$ again along
  the corresponding component of $\phi(Z)$. Along these cylinders, the
  fiberwise tangent bundle is naturally trivialized by the tangents of
  $F\cap P_t$. Along $\Pd$, the contact bundle coincides with the
  tangent bundle to $\Pd$. The result now follows.
\end{proof}

Consider next a basis subordinate to the homologically essential
knot $L\subset P$, that is,  
$\{a_1,...,a_g\}$ is a basis for $H_1(P,\partial P)$ with the property that
$a_2,...,a_g$ are disjoint from $L$ and $a_1$ meets it in a single
transverse intersection point. We can close off the arcs to get
a basis for $H_1(P)$. The Thurston--Bennequin number of $L$ can be
read off from these data as follows.

\begin{lem}
  \label{lemma:ThurstonBennequin}
  Suppose $L$ is null-homologous in $Y(\phi)$ so that, according to
  Lemma~\ref{lemma:NullHomologousKnot}, $[L]=\phi_*(Z)-Z$. Then,
  writing $Z$ in the above basis as
  $$Z=\sum_{i=1}^g n_i\cdot a_i,$$
  we find that $n_1$ is the Thurston--Bennequin invariant of $L$,
  where here we have oriented $a_1$ so that $\#L\cap a_1=+1$.
\end{lem}

\begin{proof}
  Recall that the Heegaard surface $\Sigma$ associated to the open
  book decomposition is gotten by doubling $P$ along its boundary,
  $\Sigma=P_{+1}\cup (-P_{-1})$.  We let $\alpha_i$ be the curve gotten
  by doubling $a_i$, $\alpha_i=a_i\cup {\overline {a_i}}$ and for the
  perturbed arcs $b_i$ we get $\beta_i=b_i\cup {\overline
    {\phi(b_i)}}$.
  
  Thus, the Heegaard diagram at $P_{+1}\subset \Sigma$ has a standard
  form (independent of the monodromy map). Our knot $L$ is adapted to
  the Heegaard diagram if it can be drawn on $P_{+1}$ so that it meets
  only $\alpha_1$ and $\beta_1$. The curve $L$ up to isotopy is in
  fact determined by any two points $w$ and $z$ in the two different
  components of $L-L\cap\alpha_1-L\cap\beta_1$.  $L$ can be thought of
  now as a union of two arcs, $\xi$ crossing only $\alpha_1$ and
  $\eta$ crossing only $\beta_1$.
  
  In the proof we will modify our Heegaard diagram in a way that it
  will accomodate the 0--surgery along $L$, and hence the Seifert
  surface will be visible as a periodic domain.  We stabilize $\Sigma$
  once to obtain a new Heegaard diagram which corresponds to a
  Heegaard splitting with the property that $L$ is supported inside of
  the $\beta$-handlebody. Specifically, we let $\Sigma'$ be the
  surface obtained by attaching a one--handle to $\Sigma$ with feet at
  $w$ and $z$. We introduce a new circle $\beta_0$ which is dual to
  the one--handle, and a circle $\alpha_0$ which meets $\beta_0$ at a
  single point, running through the one--handle, and completed by the
  arc $\eta$ outside the handle. The curve $\beta_0$ is a meridian for
  $L$. Let $\gamma_0$ be a circle which runs through the new
  one-handle so as not to meet $\alpha_0$, and runs along $\xi$ in the
  surface. See Figure~\ref{fig:StabOpBookHeeg} for an illustration.

  \begin{figure}[ht]
    \begin{center}
      \mbox{\vbox{\epsfbox{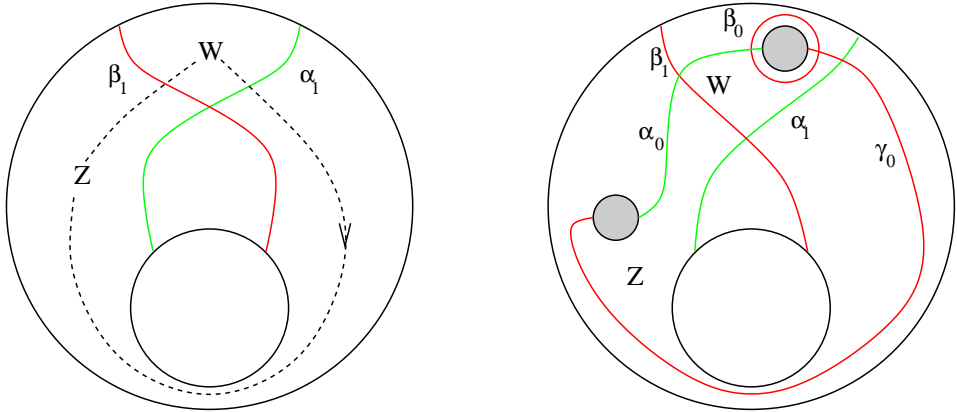}}}
    \end{center}
    \caption{
      {\bf{Stabilizing the Heegaard diagram for an open book.}}
      \label{fig:StabOpBookHeeg} On the left, we have pictured
      $P_{+1}$, thought of as an annulus, with $L$ indicated by a
      dotted line (and the two basepoints $w$ and $z$). On the right,
      we have indicated a stabilization: the two grey circles
      represent the feet of a one-handle to be added. The curves
      $\alpha_0$ and $\gamma_0$ run through this one-handle.}
  \end{figure}

  The Heegaard diagram
  $(\Sigma',\{\alpha_0,...,\alpha_g\},\{\beta_0,...,\beta_{g}\})$
  obtained by the above modification still represents $Y$. Moreover,
  $(\Sigma',\{\alpha_0,...,\alpha_g\},\{\gamma_0,\beta_1,...,\beta_g\})$
  represents the three-manifold gotten by performing surgery on $L$
  with its Thurston--Bennequin framing. Other integral surgeries on
  $L$ are represented by replacing $\gamma_0$ by a suitable smoothing
  $\delta_m$ of $\gamma_0 + m \beta_0$ with $m\in\Z$.  The
  zero--surgery is characterized by the property that there is a
  periodic domain $\Pd$, containing $\delta_m$ with multiplicity one
  along its boundary.
  
  We claim that $\delta_m$ appears as a boundary component for a
  periodic domain in $(\Sigma',\{\alpha_0,...,\alpha_g\},
  \{\delta_m,\beta_1,...,\beta_g\})$ precisely when $m=n_1$ (in the
  notation of Lemma~\ref{lemma:ThurstonBennequin}).  We construct this
  periodic domain in three pieces, $A$, $B$, and $C$ which we define
  presently.  Let ${\overline L}$ denote the copy of $L$ in
  $P_{-1}\subset \Sigma$.  The two-chain $A$ is chosen so that
  $$\partial A-\alpha_0-\gamma_0+{\overline L}\in \Span([\alpha_i]_{i=2}^g);$$ 
  the two-chain $B$ has the property that
  $$\partial B-{\overline L}\in\Span([\beta_0],
[\alpha_i-\beta_i]_{i=1}^g),$$ hence gives a relation between
${\overline L}$, and a linear combination of $\beta_0$, and the
$\{\alpha_i-\beta_i\}_{i=1}^g$. Finally, $C$ is a two--chain
connecting $\beta_0$, $\gamma_0$, and $\delta_m$, i.e.
  $$\partial C = \delta_m-\gamma_0-m\beta_0,$$
  see the picture in
  Figure~\ref{fig:RegionC}.

  \begin{figure}[ht]
    \begin{center}
      \mbox{\vbox{\epsfbox{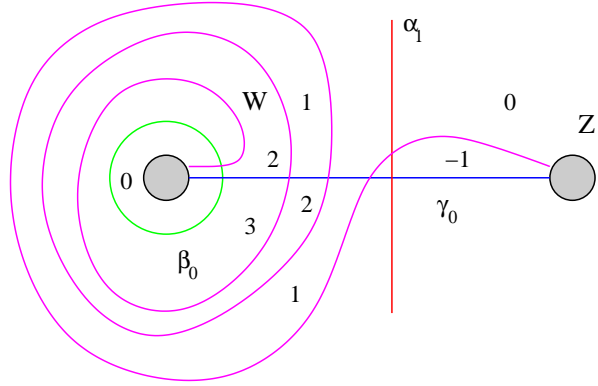}}}
    \end{center}
    \caption{
      {\bf{Two--chain of type $C$ when $m=3$.}}
      \label{fig:RegionC} Its various local multiplicities are
      indicated.}
  \end{figure}
  
  The chain $A$ exists, as follows.  Recall that $\Sigma'$ is obtained
  by stabilizing $\Sigma=P\cup (-P)$.  We see that
  $\alpha_0+\gamma_0$ is homologous in $\Sigma'$ to a copy of $L$,
  which we think of as supported in $\Sigma$ (i.e. it is supported
  away from the stabilization region in $\Sigma'$).  It suffices now
  to show that $L-{\overline L}$, thought of now as a curve in
  $\Sigma$, is homologous to a sum of the $\alpha_i$ (with $i>1$).  We
  see this as follows.  Cutting $\Page$ along $a_i$ with $i>2$, we end
  up with an annulus $X$ with some boundary arcs given by $a_i$, and
  containing $L$ as its core. Thus $L$ separates $\Page-a_2-...-a_g$
  into two components $C_1$ and $C_2$.  Similarly, if we
  cut $\Sigma=\Page\cup (-\Page )$ along the $\alpha_i=a_i\cup
  {\overline a}_i$ ($i>1$), we see that the union of $L$ and
  ${\overline L}$ separates $\Sigma-\alpha_2-...-\alpha_g$ into two
  components, $C_1\cup ( -C_1)$, and $C_2\cup (-
  C_2 )$. We let $A$ be the appropriate component (as determined by
  orientations).

  The chain $B$ is constructed from the two--chain $\Pd$ from
  Lemma~\ref{lemma:NullHomologousKnot}, and by drilling out the disk
  with multiplicity $m$ in the region between $\alpha_1$ and
  $\beta_1$, i.e. $\partial B$ contains ${\overline L}$ with
  multiplicity one, $\alpha_1$ with multiplicity $-n_1$, and so
  contains $\beta_0$ with multiplicity $n_1$.  The chain $C$ exists
  from the construction of $\delta_m$.
  
  The condition that $m= n_1$ ensures that, in the boundary of the
  sum $A+B+C$, the multiplicity of $\beta_0$ is zero. Since the
  multiplicity of $\gamma_0$ is also zero, we see that $A+B+C$
  actually represents a periodic domain for the zero--surgery, as
  claimed.
\end{proof}

\begin{proof}[Proof of Theorem~\ref{thm:SpecialCaseAlex}]
  We think of 
  $$(\Sigma,\{\delta_m,\beta_1',...,\beta_g'\},
  \{\beta_0,...,\beta_g\} \{\alpha_0,...,\alpha_g\}, z)$$ (where
  $\beta_i'$ are small isotopic translates of $\beta_i$) as a Heegaard
  triple representing the two--handle cobordism corresponding to the
  zero--surgery on $L\subset Y(\phi)$. Generators for $\Tb\cap\Tc$
  have ``nearby'' representatives in $\Td\cap\Ta$, as in
  Figure~\ref{fig:NearbyGenerator}.
  \begin{figure}[ht]
    \begin{center}
      \mbox{\vbox{\epsfbox{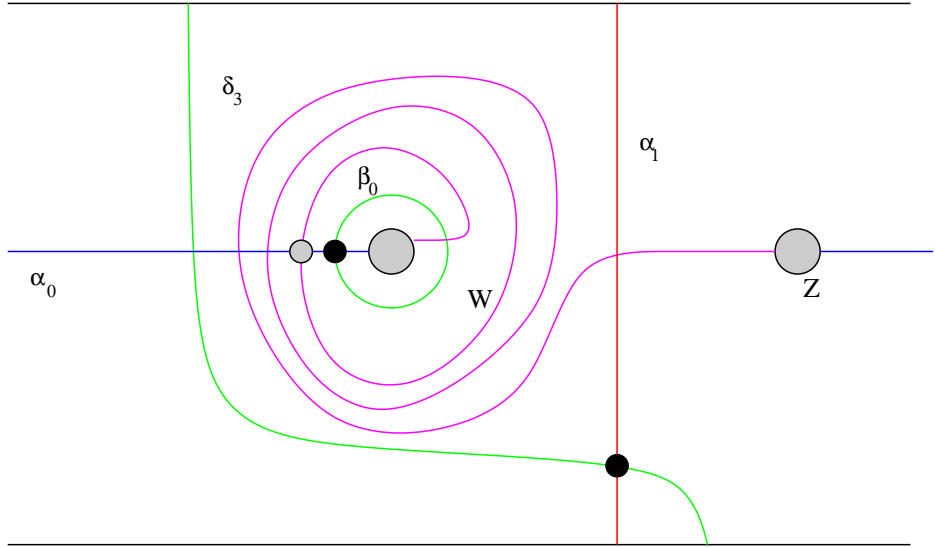}}}
    \end{center}
    \caption{
      {\bf{Nearby generator.}}
    \label{fig:NearbyGenerator}
    The large grey circles represent the handle (compare
    Figure~\ref{fig:StabOpBookHeeg}); the black circles represent the
    original generator (using $\beta_0$) while the grey one represents
    the surgery.}
  \end{figure}
  The Alexander grading of $\x$ can be calculated by the Alexander
  grading of the corresponding nearby point $\x'$, since the two can
  be connected by a triangle $\psi\in\pi_2(\Theta,\x,\x')$ with
  $n_{w}(\psi)=n_{z}(\psi)=0$ (cf.~\cite[Section~2.3]{OSzknot}).
    In turn, the Alexander grading of $\x'$ is calculated with
  the help of the formula
  $$\langle c_1(\spinc(\x')),[{\widehat F}]\rangle =
  2p_{\x'}({\widetilde \Pd})+e({\widetilde \Pd}),$$
  where here ${\widetilde \Pd}$ is a periodic domain representing the
  homology class of ${\widehat F}$,
  and  $p_{\x'}({\widetilde \Pd})$ denotes the sum of the local multiplicities
  of ${\widetilde \Pd}$ at the components of $\x'$
  cf.~\cite[Proposition~7.5]{HolDisk}.  As in the proof of
  Lemma~\ref{lemma:ThurstonBennequin}, we construct a two--chain
  representing a periodic domain in the zero--surgery, cut into three
  pieces $A$, $B$, and $C$ (in the notation from the previous proof).
  We claim that
  \begin{equation}
    \label{eq:MasA}
    2p_{\x'}(A)+e(A)=0.
  \end{equation}
  To see this, we argue as follows. Note that $p_{x_0}(A)=\OneHalf$.
  There are three types of remaining component of $x_i\in \x'$ with
  $i>0$: those which are contained in the interior of $A$, those which
  are disjoint from the closure of $A$ (like $x_1$), and those which
  are contained in the boundary of $A$. At those which are not
  contained in the closure of $A$, $p_{x_i}(A)=0$. At those which are
  in the interior, $p_{x_i}(A)=1$.  At those which are on the
  boundary, $p_{x_i}(A)=\OneHalf$. But also it is easy to see that the
  Euler measure of $A$ is given by
  $$e(A)=-1-\#\{i\big| p_{x_i}(A)=\OneHalf\}-2\#\{i\big| p_{x_i}(A)=1\};$$ 
  indeed $A$ is a connected surface of genus $\#\{i\big|
  p_{x_i}(A)=1\}$, and whose number of boundary components is given by
  $3+\#\{i>0\big|p_{x_i}(A)=\OneHalf\}$.
  Equation~\eqref{eq:MasA} now follows.

  We also claim that 
  \begin{equation}
    \label{eq:PointB}
    p_{\x'}(B)=n_1.
  \end{equation}  Indeed, $p_{x_0}(B)=n_1$ (since  $x_0$ 
  lies in the region between $\alpha_1$ and $\beta_1$), while for all
  $i>0$, $p_{x_i}(B)=0$ since at each $x_i\in\x'$ which meets $B$, the
  four corners have local multiplicities $0$, $n_i$, $0$, and $-n_i$.
  From the relationship between $B$ and $\Pd$, it is clear that
  \begin{equation}
    \label{eq:EulerB}
    e(B)=-(e(\Pd)-n_1),
  \end{equation}
  since $B$ is obtained from $\Pd$ by removing a disk with
  multiplicity $n_1$, and we evaluate on the $(-P_{-1})$--side of the
  Heegaard surface (reversing the sign of the result).

  Finally, 
  \begin{eqnarray}
    e(C)=0 &{\text{and}}&
    p_{\x'}(C)=-n_1+\OneHalf.
    \label{eq:MasC}
  \end{eqnarray}
  The first of these follows directly by inspecting the region $C$
  (see Figure~\ref{fig:RegionC}). For the second, observe that
  $p_{x_i'}(C)=0$ except when $i=0$, in which case
  $p_{x_0'}(C)=-n_{1}+\OneHalf$ (compare
  Figure~\ref{fig:NearbyGenerator}, and the domain for $C$ pictured in
  Figure~\ref{fig:RegionC}).

  Combining Equations~\eqref{eq:MasA}, \eqref{eq:PointB}, \eqref{eq:EulerB},
  and~\eqref{eq:MasC}, we conclude that
  $$\langle c_1(\spinc(\x')),[\widehat F]\rangle = e(\Pd )-n_1+1.$$ In
  view of Lemmas~\ref{lemma:RotationNumber}
  and~\ref{lemma:ThurstonBennequin}, the theorem now
  follows.
\end{proof}

Consider now the case of general ambient 3--manifold $Y$.  In this
case both the Alexander grading $A$ and the rotation number $\rot $
(as integer--valued invariants) require an additional choice: we need
to fix a Seifert surface $F$ for $L$.  Indeed, $A_F(\x )$ is defined
by the formula \eqref{e:Algra} with $\spinct = \spinct _{\x}$ being
the relative $\SpinC$ structure of $\x$, while $\rot _F (L)$ is the
integral of the relative Euler class of the oriented 2--plane bundle
$\xi $ on $F$, with the trivialization of $\xi _{\partial F}$ given by
the oriented tangent vectors of $L$. With this choice, the quantities
$A_F (\LegInv (L))$ and $\rot _F(L)$ are well--defined and we get

\begin{thm}\label{t:altalanos}
For a fixed Seifert surface $F$ of the Legendrian knot $L\subset (Y,
\xi )$ we have
\[
2A_F(\LegInv (L))=\tb (L) -\rot _F (L) +1 .
\]
\end{thm}
\begin{proof}
In Lemma~\ref{lemma:NullHomologousKnot} we can choose the
decomposition $[L]=\phi _* (Z)-Z$ in such a way that the Seifert
surface resulting from $Z$ by Lemma~\ref{lemma:RotationNumber} is
homologous to $F$ in the relative homology group $H_2 (Y,L)$. The rest
of the proof is then applies verbatim.
\end{proof}

\subsection{Maslov gradings}

Equation~\eqref{eq:MaslovGradingFormula} is much easier to establish.
In fact, we establish the following more general version:

\begin{thm}
  \label{thm:MaslovGradingFormula}
  Let $(Y,\xi)$ be a contact three--manifold with the property that
  $c_1(\xi)$ is a torsion homology class, so that
  $\CFa(Y,\spinc(\xi))$ has a rational Maslov grading, and also the
  two--plane field $\xi$ has a Hopf invariant $d_3(\xi)\in\Q$.
  Suppose moreover that $L\subset Y$ is a null--homologous Legendrian
  knot.  Then, we have that
  $$ 2A(\LegInv(L))-M(\LegInv(L))=d_3(\xi). $$
\end{thm}

For more on the absolute (rational) Maslov grading on
$\CFa(Y,\spinct)$ where $c_1(\spinct)$ is a torsion class,
see~\cite{AbsGraded}. In defininig the Hopf invariant $d_3(\xi )$ of
the 2--plane field underlying the contact structure $\xi$, we follow
the convention used in \cite{OSzcont}. Notice that when $c_1(\xi
)$ is torsion, the Alexander grading $A(\LegInv (L))$ of the
Legendrian invariant and the rotation number $\rot (L)$ of the
Legendrian knot are independent of the chosen Seifert surface.

We can continue to think of $\CFKm(Y,\spinct)$ as a bigraded group,
with a Maslov grading induced from $\CFm(Y,\spinct)$, and with an
Alexander grading (which ordinarily we think of as given by relative
$\SpinC$ structures compatible with $\spinct$) defined as half the
first Chern class of the relative $\SpinC$ structure evaluated on a
Seifert surface for $K$.  This latter quantity will be denoted by $s$.

Consider the map
$$F\colon \CFKm(-Y, L, \spinct , s )\longrightarrow \CFa(-Y, \spinct )$$
gotten by setting $U=1$, and then viewing $z$ as the basepoint for
$\CFa(-Y, \spinct )$.

\begin{prop}
  \label{prop:GradingPsi}
  The map $F$ sends
  $\CFKm_d(-Y,L,\spinct,s)\longrightarrow \CFa_{d-2s}(-Y,\spinct).$
\end{prop}

\begin{proof}
  For each given $s$, the map $F$ clearly preserves the relative
  Maslov grading, as that is given by $\Mas(\phi)$ where
  $\phi\in\pi_2(\x,\y)$ satisfies $n_{w}(\phi)=n_{z}(\phi)=0$.
  Moreover, if $\phi\in\pi_2(\x,\y)$ satisfies $n_{w}(\phi)=0$ and
  $n_{z}(\phi)=k$, so that $A(\x)-A(\y)=k$, then under the
  specialization, we see that $\phi$ drops Maslov grading by $2k$.  It
  follows at once that there is a constant $c$ with the property that
  $F$ sends $\CFKm_{d}(-Y,L,\spinct,s)$ to 
 $ \CFa_{d-2s-c}(-Y,\spinct )$.  The fact $c=0$ follows from
  conjugation symmetry of Floer homology.
\end{proof}

\begin{proof}[Proof of Theorem~\ref{thm:MaslovGradingFormula}.]
The construction of the Legendrian invariant $\LegInv (L)$ implies
that $F(\LegInv(L))=c(Y, \xi)$, where here $c(Y, \xi)\in\CFa(-Y)$ is
the contact invariant from~\cite{OSzcont}. By \cite{OSzcont} the
Maslov grading of the contact invariant is $M(c(\xi))=-d_3(\xi)$,
hence Proposition~\ref{prop:GradingPsi} above implies that
$M(\LegInv(L))-2A(\LegInv(L))=-d_3(\xi)$, concluding the proof.
\end{proof}
\begin{proof}[Proof of Theorem~\ref{thm:BigradingFormula}]
  Theorem~\ref{thm:BigradingFormula} is now a combination of
  Theorems~\ref{thm:SpecialCaseAlex}
  and~\ref{thm:MaslovGradingFormula}.
\end{proof}

\section{Transverse simplicity of knots in $(S^3, \xi _{st})$}
\label{s:ecknots}

\subsection{The Eliashberg--Chekanov knots}

\newcommand{\cva}{\underline {{\bf {\rm c}}}}

\begin{figure}[ht] 
\begin{center} 
\setlength{\unitlength}{1mm}
{\includegraphics[height=4cm]{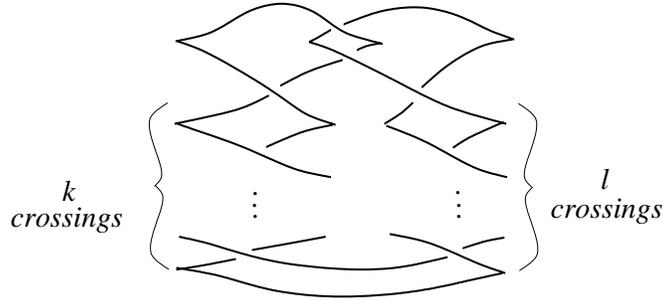}}
\end{center} 
\caption{{\bf The Legendrian knots $E(k,l)$.} These knots are smoothly
  isotopic to $E_n$, with $k+l-1=n$.}
\label{f:chekan} 
\end{figure}

We will demonstrate the power of the transformation rule  proved in 
Section~\ref{s:surg} by computing the invariants of the Eliashberg--Chekanov
Legendrian knots and verify Theorem~\ref{t:ecnonsimple}.  To this end,
consider the Legendrian knot $E(k,l)$ given by Figure~\ref{f:chekan}.

\begin{prop}\label{p:szamok}
The knot $E(k,l)$ is a Legendrian knot in the standard contact
3--sphere smoothly isotopic to the Eliashberg--Cheknov knot $E_n$ with
$k+l-1=n$, cf. \cite{EFM}. The Thurston--Bennequin and rotation
numbers of $E(k,l)$ are given as $\tb(E(k,l))=1$ and $\rot
(E(k,l))=0.$ \qed
\end{prop}

\begin{cor}
The Legendrian invariant $\LegInv (E(k,l))$ is a nonzero element
of the knot Floer homology group $\HFKm _{2}(-S^3, E_n, 1)$.
\end{cor}
\begin{proof}
According to the nonvanishing of the invariant in a Stein fillable contact
3--manifold, we get $\LegInv (E(k,l))\neq 0$. The Alexander and Maslov
gradings of $\LegInv (E(k,l))$ can be computed from the rotation and 
Thurston--Bennequin numbers (given in Proposition~\ref{p:szamok}) through
the formulae of Theorem~\ref{thm:BigradingFormula}.
\end{proof}

Let $k=\frac{n+1}{2}$.
The knot Floer homology group $\HFKm (-S^3, E_n)=\HFKm (S^3, m(E_n))
=\HFKm (m(E_n))$ for $n\geq 1$ and odd is given as
$$\HFKm_M (m(E_n), A)=\left\{\begin{array}{lll}
\Field ^{k} & (A=1, M=2)\\
\Field ^{k} & (A=0, M=1)\\
\Field    &   (A=-i\leq 0, M=-2i).
\end{array}\right.$$
Also, $\HFKa (m(E_n))$ (which can be read off directly 
from the Alexander polynomial and the signature of the knot since
$E_n$ is an alternating knot, see~\cite{AltKnot,2Bridge}) is given as 
$$\HFKa_M (m(E_n), A)=\left\{\begin{array}{lll}
\Field ^k & (A=1, M=2)\\
\Field ^n & (A=0, M=1)\\
\Field^k    &   (A=-1, M=0).
\end{array}\right.$$

\begin{cor} The Legendrian invariant $\LegInva (E(k,l))$ is a nonzero
element of the knot Floer group $\HFKa (m(E_n))$.
\end{cor}
\begin{proof}
The explicit description of the specialization map $\HFKm (m(E_n))
\to \HFKa (m(E_n))$ (when setting $U=0$) and the fact that the
invariant lives in bidegree $(A=1, M=2)$ readily implies the
corollary.
\end{proof}

We would like to verify that the Legendrian knots $E(k,l)$ with
$k+l-1=n$, $k,l$ odd and $k\leq l$ have different Legendrian
invariants.  As usual, we use  $\HFKa$, implying the
similar distinction result for the invariants in the $\HFKm$--groups.

\begin{thm}\label{t:chekk} 
Let us fix a knot $E_n \subset S^3$.  There are identifications of
$\HFKa (S^3, m(E(k,l)))$ with $\HFKa (S^3, m (E_n))$, for $k+l-1=n$,
$k,l$ odd such that the images of the the Legendrian invariants
$\LegInva (E(k,l))$ and $\LegInva (E(k',l'))$ are equal in $\HFKa
(S^3, m(E_n))$ if and only if $k=k'$ and $l=l'$.
\end{thm}
The action of the mapping class group taken into account, this
statement gives
\begin{proof}[Proof of Theorem~\ref{t:ecnonsimple}]
Since $E_n$ is a two--bridge knot, which is not a torus knot, and by
\cite{HT} all non--torus 2--bridge knots are hyperbolic,
\cite[Theorem~2.7]{RW} implies that the mapping class group $\MCG
(S^3, E_n)$ is isomorphic to $\Z /2\Z$.  The Legendrian knots $E(k,l)$
and $E(l,k)$ are Legendrian isotopic \cite[Theorem~4.1]{EFM}, hence
under the action of the mapping class group the element $\LegInva
(E(k,l))$ is mapped to $\LegInva (E(l,k))$. This identifies the
$\Z/2\Z$--action of $\MCG (S^3, E_n)$ on the knot Floer homology group
$\HFKa _2(m(E_n),1)$, and shows that $\LegInva (E(k,l))$ and $\LegInva
(E(k',l'))$ are in different $\MCG (Y,L)$--orbits, provided
$k,l,k',l'$ are odd, $k+l-1=k'+l'-1=n$, $k\leq l, k'\leq l'$ and
$k\neq k'$. This shows that the corresponding Legendrian knots and
their negative stabilizations are not isotopic, concluding the proof
of the theorem.
\end{proof}
\begin{rem} {\rm Notice that the above result, in conjunction with
    \cite[Theorems~2.2~and~4.2]{EFM} shows that among the knots given
    by positive transverse push--offs of the diagram of
    Figure~\ref{f:chekan} there are exactly $\lceil \frac{n}{4}\rceil$
    transversely non--isotopic. It is still an open question whether
    there are further transverse representatives of $E_n$ with
    self--linking number $1$ not transverse isotopic to any of the
    transverse push--offs of the Legendrian knots considered above.}
\end{rem}

\begin{figure}[ht] 
\begin{center} 
\setlength{\unitlength}{1mm}
{\includegraphics[height=5cm]{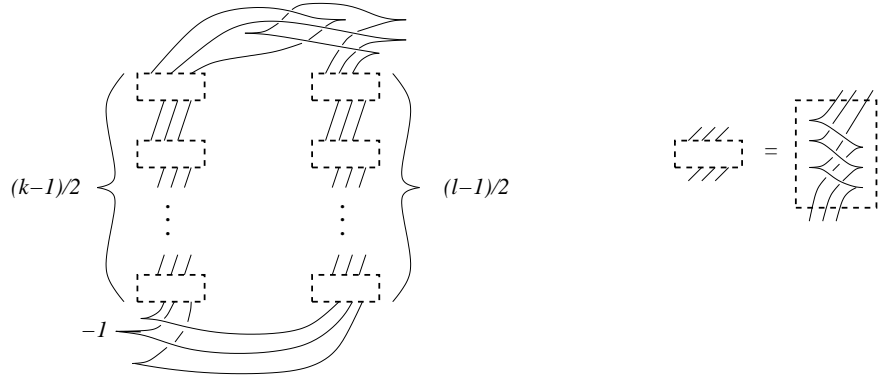}}
\end{center} 
\caption{{\bf The Legendrian knot $E(k,l)$ together with an unknot.}}
\label{f:surgery} 
\end{figure}

For the proof of Theorem~\ref{t:chekk}, consider the 2--component
Legendrian link of Figure~\ref{f:surgery}. Notice that the linking 
number of the two knots is zero.
\begin{figure}[ht] 
\begin{center} 
\setlength{\unitlength}{1mm}
\includegraphics{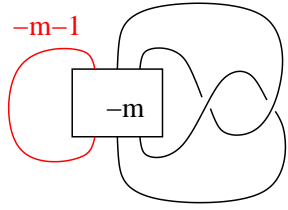}
\end{center} 
\caption{{\bf The smooth knots of Figure~\ref{f:surgery}.} 
We let $m=\frac{n+1}{2}$.}
\label{f:smoothsurgery} 
\end{figure} 
The smooth diagram underlying
Figure~\ref{f:surgery} is given by Figure~\ref{f:smoothsurgery}.
Applying contact $(-1)$--surgery on the unknot component of
Figure~\ref{f:surgery}, we get a Legendrian knot $E'(k,l)$ in the lens
space $L(m+1,1)$ with contact structure $\xi _{k,l}$. Let $S$ denote
the Legendrian push--off of the unknot in Figure~\ref{f:surgery}. It
is a standard fact in contact topology that contact $(+1)$--surgery
along $S$ cancels the first contact $(-1)$--surgery, and hence
provides the standard contact 3--sphere with the Legendrian knot
$E(k,l)$ in it. The surgery along $S$ induces the map
\[
\Fa _{S}\colon \HFKa (-L(m+1, 1), E'(k,l))\to \HFKa (-S^3, E(k,l)) .
\]
(Notice the orientation reversal on the 3--manifolds.)

\begin{prop} \label{p:dist}
The Legendrian invariants $\LegInva(E'(k,l))$ are all distinct, 
and the map $\Fa _{S}$ is injective on the subgroup of the knot
Floer homology with Alexander grading 1.
\end{prop}
\begin{proof}
The surgery along $S$, with orientation reversed, gives rise to a
distinguished triangle of knots, depicted in
Figure~\ref{f:triangle}. For simplicity, here we slid the surgery
curve $S$ over the unknot to which it is a Legendrian push--off,
resulting in a meridional unknot with framing 0.
\begin{figure}[ht] 
\begin{center} 
\setlength{\unitlength}{1mm}
{\includegraphics[height=8cm]{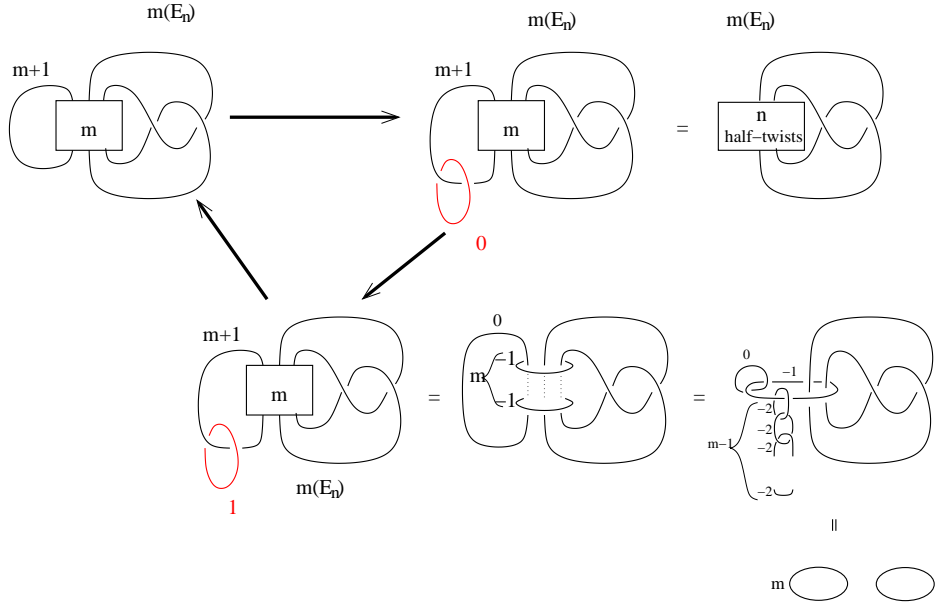}}
\end{center} 
\caption{{\bf The distinguished triangle of knots induced by the surgery
along $S$ (after the reversal of orientation).} Recall that
$m=\frac{n+1}{2}$.}
\label{f:triangle} 
\end{figure} 
As it was explained in Subsection~\ref{ss:mapsbysurgery}, such a
distinguished triangle of knots induces an exact triangle on the knot
Floer homology groups. In addition, since the surgery curve and the
knot under inspection gives zero linking, we can consider groups only
with a fixed Alexander grading, since all maps do respect that
Alexander gradings. Since $\LegInva (L)$ is a nonzero element with
Alexander grading $1$, we examine the $A=1$ groups only.

The third term of the triangle is an unknot in the lens space
$-L(m,1)$, therefore at Alexander grading 1 the corresponding knot
Floer group vanishes, implying that the map $\Fa _{S}$ is an
isomorphism on that particular Alexander grading.

The element $\LegInv (E'(k,l))$ specializes to $c(L(m+1,1), \xi
_{k,l})$ under the specialization $U=1$. Since $\xi _{k,l}$ and $\xi
_{k',l'}$ induce the same $\SpinC$ structure only if $k=k'$ and
$l=l'$, we conclude that the invariants $\LegInv (E'(k,l))$ are
different for different $k$'s.
\end{proof}

\begin{proof}[Proof of Theorem~\ref{t:chekk}]
The injectivity of $\Fa _{S}$ and the fact that all $\LegInva (E' (k,l))$
are different, together with the naturality formula
\[
\Fa _{S}(\LegInva (E'(k,l))=\LegInva (E(k,l))
\]
concludes the proof.
\end{proof}

\subsection{Two--bridge knots}

The method of the above argument can be generalized to find
further examples of knot types containing distinct transverse knots
with identical self--linking numbers.  Here we formulate a
result exploiting the same ideas used above, and then provide some
further families where this principle can be used.  The candidates
will be constructed through the Legendrian satellite construction
described in \cite{Lenny} --- in fact, the Eliashberg--Chekanov knots
considered above are special cases of this construction.

Let us start by recalling the Legendrian satellite construction.  To
this end, let ${\tilde {L}}$ denote a Legendrian link in $S^1\times
\bfr ^2$, which can be conveniently depicted by cutting $S^1$ open at
a point, hence ${\tilde {L}}$ can be drawn in a box,
cf. \cite[Figure~22]{Lenny}.  Now for a Legendrian knot $L$ consider
its standard neighbourhood.  By an appropriate contactomorphism
between this solid torus and the one containing ${\tilde {L}}$ we can
embed ${\tilde {L}}$ into the neighbourhood of $L$. We define
$S(L,{\tilde {L}})$ as this new Legendrian knot. If $w({\tilde {L}})$
denotes the winding number of ${\tilde {L}}$ in $S^1\times \bfr ^2$,
then we have
\begin{equation}
\label{eq:TBeq}
\tb (S(L, {\tilde {L}}))=(w({\tilde {L}}))^2\tb (L)+\tb ({\tilde {L}})
\end{equation}
and
\begin{equation}
\label{eq:Masik}
\rot (S(L, {\tilde {L}}))=(w({\tilde {L}}))^2\rot (L)+\rot ({\tilde {L}}).
\end{equation}
Consequently, in case $w({\tilde {L}})=0$, the Thurston--Bennequin
and the rotation numbers of $S(L,{\tilde {L}})$ are independent
of $L$.

\begin{figure}[ht] 
\begin{center} 
\setlength{\unitlength}{1mm}
{\includegraphics[height=4cm]{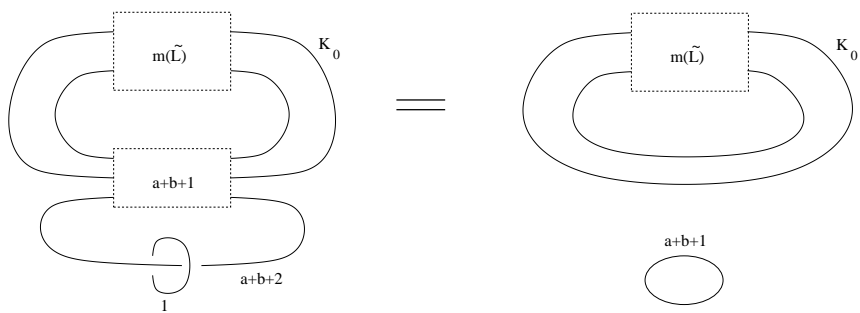}}
\end{center} 
\caption{{\bf The third knot $K_0$ in the distinguished triangle of knots.}}
\label{f:alt}
\end{figure}

Suppose that ${\tilde {L}}$ is a Legendrian solid--torus knot with
$w({\tilde {L}})=0$, and $U(a,b)$ is a Legendrian realization of the
unknot, with $\tb(U(k,l))=-1-(a+b)$ and $\rot(U(a,b))=a-b$. Let
$L_{k,l}=S(U(a,b), {\tilde {L}})$ denote the Legendrian satellite of
$U(a,b)$ with ${\tilde {L}}$, where here $k=2a+1$ and $l=2b+1$.  Let
$K_0$ denote the knot (in some lens space) given by the sugery
description of Figure~\ref{f:alt}.
Let $g(K_0)$ denote genus of the knot $K_0$. 
(By~\cite{YiNi}, this quantity bounds the largest Alexander grading
with nontrivial knot Floer homology in $\HFKm (Y, K_0)$.)

\begin{thm}\label{t:alt}
Let $k=2a+1$ and $l=2b+1$. If $\tb ({\tilde {L}})+\rot ({\tilde {L}})
>2g(K_0)-1$ and the symmetry group of the smooth knot underlying
$L_{k,l}\subset S^3$ is of order $t$ then the knot type of $L_{k,l}$
admits at least $ \lceil\frac{k+l-1}{2t}\rceil$ transversely
non--isotopic transverse representatives.
\end{thm}

\begin{proof}
  The same set--up as for the Eliashberg--Chekanov knots provides
  the Legendrian knots $L'_{k,l}$ in the lens space $L(a+b+2,1)$, and
  the map induced by the $(+1)$--surgery is an isomorphism again,
  since in the surgery triangle the third term vanishes. This
  vanishing is because in the Alexander grading
  $A=\frac{1}{2}(\tb(L_{k,l})+\rot(L_{k,l})+1)= \frac{1}{2}(\tb
  ({\tilde {L}})+\rot ({\tilde {L}})+1)$ the knot Floer homology group
  of $K_0$ vanishes by the assumption. Now the analogue of
  Proposition~\ref{p:dist} provides different invariants before the
  action of the mapping class group $\MCG(S^3, L_{k,l})$ is taken into
  account. Since in $\HFKa(S^3, m (L_{k,l}))/\MCG (S^3,L_{k,l})$ we
  will still have at least $\lceil \frac{k+l-1}{2t}\rceil$ different
  invariants, the proof follows.
\end{proof}

\begin{figure}[ht] 
\begin{center} 
\setlength{\unitlength}{1mm}
{\includegraphics[height=4cm]{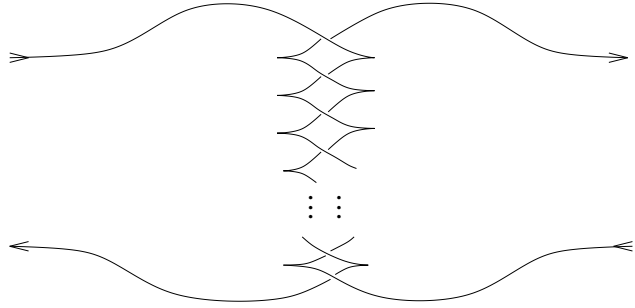}}
\end{center} 
\caption{{\bf Further examples.}}
\label{f:tovabbi}
\end{figure}

A simple example of ${\tilde {L}}$ with $w({\tilde {L}})=0, \tb=1$ and
$K_0$=Unknot is given by Figure~\ref{f:tovabbi}. Notice that the
orientation depicted in the figure implies that there are an even
number of crossings in the projection.
Since the knots $S(U(a,b), {\tilde {L}})$ are all 2--bridge
knots, and those knots have small symmetry groups provided they are
not torus knots, we find further examples of transversely non--simple
families of knot types.  

\begin{figure}[ht] 
\begin{center} 
\setlength{\unitlength}{1mm}
{\includegraphics[height=4cm]{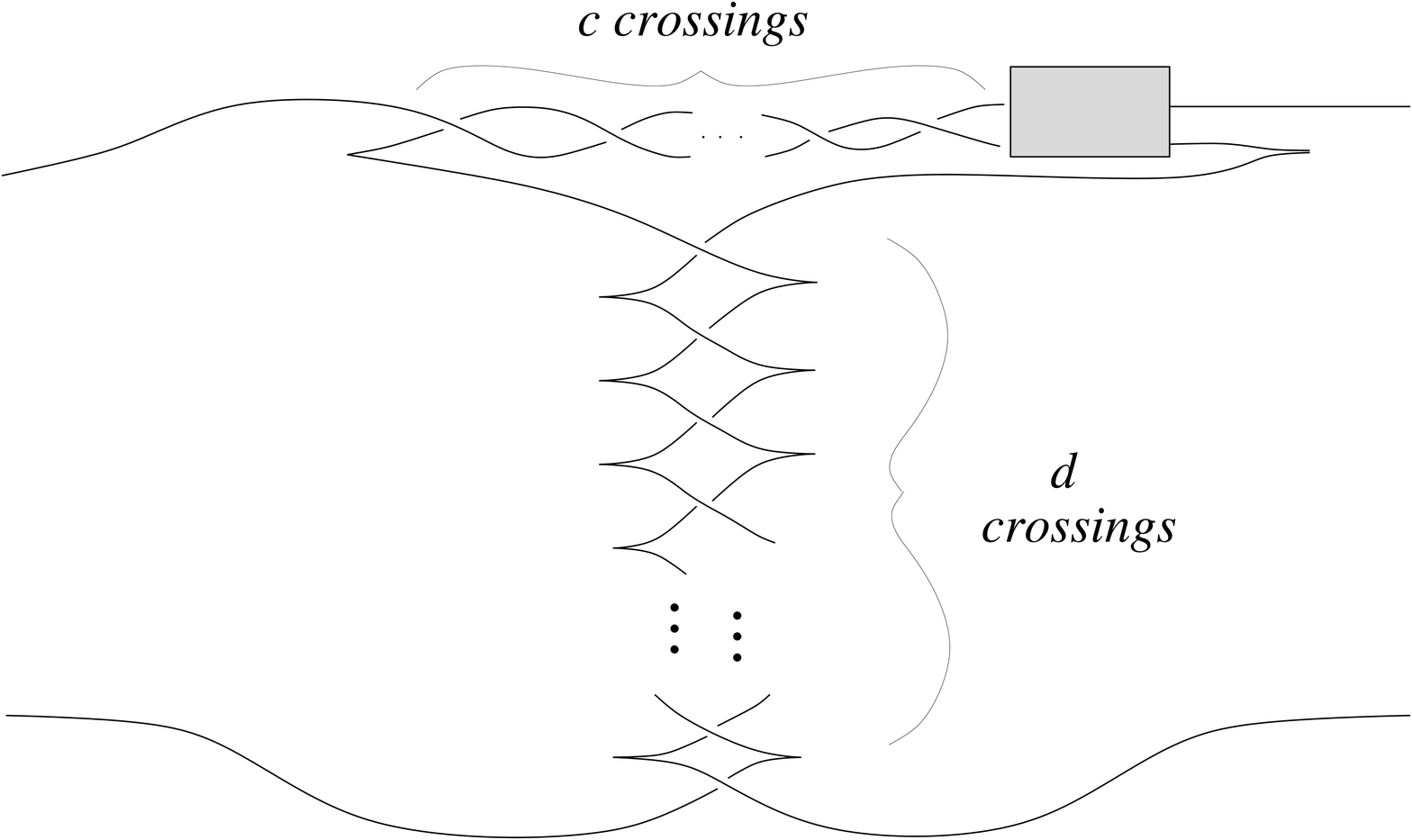}}
\end{center} 
\caption{Iterating the construction}
\label{f:modif}
\end{figure}

There are various ways to further generalize this construction. For
example, the top crossing of Figure~\ref{f:tovabbi} can be replaced
with a sequence of crossings, as it is shown in Figure~\ref{f:modif}.
In order to get a knot, rather than a link, we require the parity of
the number $c$ of these new crossings to be odd. The
Thurston--Bennequin number of the knot $S(U(a,b), {\tilde{L}})$ with
${\tilde {L}}$ given by Figure~\ref{f:modif} (with strands simply
passing through the gray box) can be easily computed to be equal to
$c$, while the Euler characteristic of a Seifert surface of $K_0$
appearing in Theorem~\ref{t:alt} is $1-c$. (Notice that $K_0$ is, in
fact, isotopic to the $(2,c)$ torus knot.) Since $L_{k,l}=S(U(a,b),
{\tilde{L}})$ with $k=2a+1$ and $l=2b+1$ for $k+l>2$ is a 2--bridge
knot which is not a torus knot, its mapping class group is known to be
isomorphic to $\bfz/2\bfz$, hence Theorem~\ref{t:alt} shows that the
knot types appearing in this construction (with $k+l>2$, both odd) are
not transversely simple.  The above construction admits a further
generalization as follows:

\begin{thm}
  Suppose that $\frac{p}{q}\in \Q$ has the continued fractions expansion
  $$\frac{p}{q} = [a_1,...,a_{2m+1}]
        = a_1 + \frac{1}{a_2 + \frac{1}{\ddots +\frac{1}{a_{2m+1}}}},
  $$
  and suppose that
  \begin{itemize}
  \item $a_2$ and $a_3$ are odd, 
  \item all $a_i$ for $i\neq 2,3$ are even and
  \item $m\geq 1$.
  \end{itemize}
  Let $K=K(p,q)$ be the corresponding two-bridge knot.
  Then, $K$ admits at least $\lceil \frac{a_1}{4}\rceil$ 
  distinct transverse realizations with self--linking number $\sln$ equal to
  $$\sum_{i=1}^m a_{2i+1}.$$
\end{thm}

\begin{proof}

  Consider the sequence ${\widetilde L}_i$ ($i=1,\ldots ,m$) of solid
  torus knots given by Figure~\ref{f:modif}, where the parameters
  $c_i,d_i$ are chosen so that $c_1,d_1$ are both odd, while $c_i,d_i$
  are even for $i>1$. Starting with $i=m$, copy ${\widetilde L}_{i}$
  into the gray box of ${\widetilde L}_{i-1}$. The parity assumptions
  on $c_i$, $d_i$ ensure that ${\widetilde L}_i$ can be drawn in the
  indicated box of ${\widetilde L}_{i-1}$ with consistent
  orientations.  At the end of this process we get a solid torus link
  ${\widetilde L}$, with $w({\widetilde L})=0$.  We use it to define
  the set of Legendrian knots $L_{k,l}=S(U(a,b),{\widetilde L})$ with
  $k=2a+1$, $l=2b+1$ as before.

  To calculate the Thurston-Bennequin number, we argue as follows.
  According to Equation~\eqref{eq:TBeq}, $\tb(S(L,{\widetilde
    L}))=\tb({\widetilde L})$.  To calculate this, observe first that
  all crossings in the diagram are positive; hence the writhe is the
  total number of crossings. Moreover, in the ${\widetilde L}_i$
  template, the $d_i$ crossings cancel with the left cusps, leaving
  only the $c_i$ to contribute to the Thurston-Bennequin invariant.
  Consequently, 
  $$\tb(L_{k,l})=\sum_{i=1}^m c_i.$$

  Applying Seifert's algorithm, we see that the negative of the Euler
  characteristics of a Seifert surface for the knot $K_0$ appearing in
  Theorem~\ref{t:alt} is equal to $\left(\sum _{i=1}^m c_i\right)-1$.
  The rotation number $\rot({\widetilde L})$ is zero, therefore by
  Equation~\eqref{eq:Masik}, we see that $\rot(L_{k,l})=0$ as well.
  This implies that the self--linking number $\sln (T_{k,l})$ of the
  transverse push--off of $L_{k,l}$ is equal to $\sum_{i=1}^m c_i.$

  In fact, it is easy to see that $L_{k,l}$ is the two-bridge knot
  $K(p,q)$ with $\frac{p}{q}=[a_1,...,a_{2m+1}]$, where $a_1=k+l$,
  $a_{2i}=d_i$, $a_{2i+1}=c_i$. Since $m\geq 1$, this is not a torus
  knot. Thus, applying~\cite{HT} and~\cite{RW} as before, we conclude
  that the mapping class group has order two. Hence,
  Theorem~\ref{t:alt} provides the stated result.
\end{proof}

\section{Mapping class group actions}
\label{sec:MCG}

We construct here the mapping class group action on knot Floer
homology.  Our discussion here is built on the constructions
from~\cite{HolDiskFour} (which dealt, however, with Heegaard Floer
homology for closed three-manifolds), combined with~\cite{OSzknot}.

A {\em marked Heegaard diagram} for a pointed knot $(Y,K,p)$ is a Heegaard
diagram 
$$(\Sigma,\{\alpha_1,...,\alpha_g\},\{\beta_1,...,\beta_{g}\},w,z)$$ for $Y$
so that $w$ and $z$ determine $K$, and $w$ corresponds to the marked point $p\in K$. We can associate a
Heegaard Floer complex $\CFKm(\Sigma,\alphas,\betas,w,z)$ to this doubly-pointed Heegaard diagram, provided that it satisfies a suitable
weak admissibility condition, see~\cite{HolDisk}.

\begin{prop}
  \label{prop:HeegaardMoves} Suppose that
  $(\Sigma^1,\alphas^1,\betas^1,w^1,z^1)$ and
  $(\Sigma^2,\alphas^2,\betas^2,w^2,z^2)$ are two marked
  Heegaard diagrams for $(Y,K)$. Then, these two diagrams can be
  connected by Heegaard moves in the sense that there are the
  following data: \begin{itemize} \item marked diagrams
    $(\Sigma,{\alphas^3},{\betas^3},w,z)$ and
    $(\Sigma,{\alphas^4},{\betas^4},w,z)$ obtained by
    stabilizing $(\Sigma^1,\alphas^1,\betas^1,w^1,z^1)$ and
    $(\Sigma^2,\alphas^2,\betas^2,w^2,z^2)$ respectively;
  \item a sequence of handleslides and isotopies from
    $(\Sigma,\alphas^3,\betas^3,w,z)$ to
    $(\Sigma,\alphas^4,\betas^4,w,z)$ respectively, so
    that none of the $\alpha$-curves crosses either $w$ or $z$
    during the handleslides and isotopies.
  \end{itemize}
\end{prop}

\begin{proof}
	This is a modification of the usual 
	Reidemeister-Singer theorem, stating that two Heegaard
        diagrams for $Y$ can be connected by stabilizations,
        de--stabilizations, handleslides, and isotopies.
	We fix a Morse function $f_0$ near $K$, with one index $3$ and
	one index $0$ critical point on $K$, and consider extensions
	of this fixed Morse function to $Y$.
\end{proof}

Consider the following definition (the {\em map defined by a strong
equivalence} in the sense of~\cite[Lemma~2.13]{HolDisk}).

\begin{defn}
  Suppose that $(\Sigma,\alphas^3,\betas^3,w,z)$ and
  $(\Sigma,\alphas^4,\betas^4,w,z)$ are admissible Heegaard
  diagrams for $Y$ which differ only by handleslides and isotopies. 
  Let $\alphas^5$,
  ${\betas^5}$ be isotopic translates of the $\alphas^3$, $\betas^3$,
  so that $(\Sigma,\alphas^5,\betas^5,w,z)$ is also an admissible diagram for $Y$,
  and   $(\Sigma,\alphas^4,\alphas^5,w,z)$ and
  $(\Sigma,\betas^5,\betas^4,w,z)$ are both weakly admissible Heegaard
  diagrams for $\#^{g}(S^2\times S^1)$.  We define a map,
  up to overall multiplication by $\pm 1$,
  $$\Phi_{3,4}\colon \CFKm(\Sigma,\alphas^3,\betas^3,w,z)
  \longrightarrow\CFKm(\Sigma,\alphas^4,\betas^4,w,z) $$
  as a composite
  $$
  \begin{CD} \CFKm(\Sigma,\alphas^3,\betas^3,w,z) @>{\Gamma}>>
    \CFKm(\Sigma,\alphas^5,\betas^5,w,z)
    @>{\Theta_{\alphas^4\alphas^5}\otimes\cdot\otimes
      \Theta_{\betas^5\betas^4}}>>
    \CFKm(\Sigma,\alphas^4,\betas^4,w,z), \end{CD} $$
  where here the
  first map is induced by isotopies, and defined using continuation
  maps (i.e. counting holomorphic disks with time-dependent boundary
  conditions), while the second is defined by counting
  pseudo-holomorphic trianlges, cf.~\cite[Section~2.3]{HolDiskFour}.
\end{defn}

Given any two marked Heegaard diagrams
$(\Sigma^1,\alphas^1,\betas^1,w^1,z^1)$ and
$(\Sigma^2,\alphas^2,\betas^2,w^2,z^2)$, for suitable choices
of almost-complex structure, the stabilization/destabilization maps
give identifications \begin{align*} f_{1,3}&\colon
  \CFKm(\Sigma^1,\alphas^1,\betas^1,w^1,z^1) \longrightarrow
  \CFKm(\Sigma,\alphas^3,\betas^3,w,z)
  \\ f_{4,2}&\colon\CFKm(\Sigma,\alphas^4,\betas^4,w,z)
  \longrightarrow
  \CFKm(\Sigma^2,\alphas^2,\betas^2,w^2,z^2).  \end{align*} Define the map
$\Psi _{1,2}$ by $f_{4,2}\circ\Phi _{3,4}\circ f_{1,3}$. According to
the following variant of~\cite[Theorem~2.1]{HolDisk}, the induced map
\[
(\Psi _{1,2})_*\colon
\HFKm(\Sigma^1,\alphas^1,\betas^1,w^1,z^1)\longrightarrow
\HFKm(\Sigma^2,\alphas^2,\betas^2,w^2,z^2)
\]
is independent of
all the choices made; that is

\begin{thm}(\cite[Theorem~2.1]{HolDisk})
  If $(\Sigma^1,\alphas^1,\betas^1,w^1,z^1)$ and
  $(\Sigma^2,\alphas^2,\betas^2,w^2,z^2)$ are two marked
  Heegaard diagrams for $(Y,K)$, then the connecting Heegaard moves from
  Proposition~\ref{prop:HeegaardMoves} induce a chain map
  $$\Psi_{1,2}\colon
  \CFKm(\Sigma^1,\alphas^1,\betas^1,w^1,z^1) \longrightarrow
  \CFKm(\Sigma^2,\alphas^2,\betas^2,w^2,z^2) $$ whose chain
  homotopy class, up to multiplication by $\pm 1$, is independent of the intermediate stages. \qed
\end{thm}

\begin{defn}
  A homeomorphism of Heegaard diagrams is a map
  $$f\colon (\Sigma^1,\alphas^1,\betas^1,w^1,z^1)\longrightarrow
  (\Sigma^2,\alphas^2,\betas^2,w^2,z^2)$$
  which is a homeomorphism
  from $\Sigma^1$ to $\Sigma^2$, carrying the set
  $\alphas^1=\{\alpha^1_j\}_{j=1}^{g}$ to the set
  $\alphas^2=\{\alpha^2_j\}_{j=1}^{g}$, $\betas ^1$ to $\betas ^2$, 
  $w^1$ to $w^2$, and $z^1$ to $z^2$.
\end{defn}

A homeomorphism 
$$f\colon (\Sigma^1,\alphas^1,\betas^1,w^1,z^1)
\longrightarrow
(\Sigma^2,\alphas^2,\betas^2,w^2,z^2)$$
of Heegaard diagrams induces a continuous map
from $\Sym^g(\Sigma^1)$ to $\Sym^g(\Sigma^2)$ carrying ${\mathbb T}_{\alpha^1}$ and ${\mathbb T}_{\beta^1}$ to
${\mathbb T}_{\alpha^2}$
and ${\mathbb T}_{\beta^2}$
respectively.  This induces a map of chain complexes
$$\Xi_{1,2}(f)\colon \CFKm(\Sigma^1,\alphas^1,\betas^1,w^1,z^1)
\longrightarrow
\CFKm(\Sigma^2,\alphas^2,\betas^2,w^2,z^2).
$$

\begin{defn}
	\label{def:MapOfHomeo}
        Suppose that $F\colon (Y,K) \longrightarrow (Y,K)$ is a
        homeomorphism (fixing $K$ pointwise). Choose a Heegaard diagram
        $(\Sigma^1,\alphas^1,\betas^1,w^1,z^1)$ for $(Y,K,p)$. The map
        $F$ induces another homeomorphic diagram
\[
(\Sigma^2,\alphas^2,\betas^2,w^2,z^2)=(F(\Sigma ^1), F(\alphas ^1),
F(\betas ^1), F(w^1),F(z^1)),
\]
together with an induced homeomorphism $$f\colon
(\Sigma^1,\alphas^1,\betas^1,w^1,z^1) \longrightarrow
(\Sigma^2,\alphas^2,\betas^2,w^2,z^2).$$ We define the action of $F$ on
the (projectivized) Floer homology of $(Y,K)$ to be the composite
        $$ \begin{CD} \HFKm(\Sigma^1,\alphas^1,\betas^1,w^1,z^1)
          @>{\Xi_{1,2}(f)_*}>> \HFKm(\Sigma^2,\alphas^2,\betas^2,w^2,z^2)
          @>{(\Psi_{2,1})_*}>>
          \HFKm(\Sigma^1,\alphas^1,\betas^1,w^1,z^1). \end{CD} $$ The
        resulting map will be denoted by
\[
A(F)\colon \HFKm(Y,K)/\pm 1 \longrightarrow \HFKm(Y,K)/\pm 1.
\]
\end{defn}

\begin{thm}
	\label{thm:Naturality}
        The above map
\[
F\in\DiffP(Y,K)\mapsto
        A(F)\in\Aut(\HFKm(Y,K)/\pm 1)
\]
 descends to a well-defined action of the mapping class group of
 $(Y,K)$ on knot Floer homology.
\end{thm}

\begin{proof}   
        It is clear from the construction that if $F_1,F_2\in\DiffP(Y,K)$, then
        $A(F_1)\circ A(F_2)=A(F_1\circ F_2)$. 

        Next, we claim that if $F_t\colon (Y,K)\longrightarrow (Y,K)$ 
        is an isotopy with  the property that $F_0$ is the identity
        map, then there is some $\epsilon>0$ with the property that
        for all $|t|<\epsilon$, $A(F_t)$ acts as the identity map.
        
        We argue as follows. Choose an allowed generic almost-complex
        structure, so that for any $\x,\y\in\Ta\cap\Tb$, and any
        non--constant $\phi\in\pi_2(\x,\y)$ with the property
        $\Mas(\phi)=0$, the moduli space $\ModFlow(\phi)$ is
        empty.  It follows that in the continuation map $\Gamma$
        (defined by counting time-dependent boundary conditions), if
        $t$ is sufficiently small and the moduli space for the
        continuation $\ModFlow(\x,\y')$ (for $\x\in\Ta\cap\Tb$ and
        $\y'\in\Ta'\cap\Tb'$) is non--empty, then the homotopy class
        $\phi$ must correspond to the constant map.  It follows that
        $\Gamma\colon
        \CFKm(\Sigma,\alphas^1,\betas^1,w,z)\longrightarrow
        \CFKm(\Sigma,\alphas^2,\betas^2,w,z)$ is the closest-point
        map, provided that $\alphas^2$ and $\betas^2$ are sufficiently
        close to $\alphas^1$ and $\betas^1$; and also that
        $\ModFlow(\x',\y')$ is empty for all $\phi\in\pi_2(\x',\y')$
        with $\Mas(\phi)=0$, where $\x',\y'\in{\mathbb
          T}_{\alpha^2}\cap {\mathbb T}_{\beta^2}$.
        
        Now, we choose $\alphas^3,\betas^3$ to be an exact Hamiltonian
        translate of the pair $\alphas^2, \betas^2$. We claim that in
        this case, for sufficiently small translates, the map induced
        by triangles is also a nearest point map. This is seen by
        identifying the map from the continuation map $\Gamma$
        (associated to the isotopy of $\alphas^3$ with $\alphas^2$ and
        $\betas^3$ with $\betas^2$) with the map defined by counting
        pseudo-holomorphic triangles, and then appealing to the
        previous paragraph. The identification of $\Gamma$ and the
        triangle map (in a related context) can be found
        in~\cite[Proposition 11.3]{LipshitzCyl}; we sketch this
        argument here. The continuation map (in the case where only
        the $\beta$--circles are moving, while the $\alpha$--circles
        stay fixed, for notational simplicity) can be thought of as
        counting pseudo-holomorphic triangles with one corner smoothed
        out, which map one boundary edge to $\Ta$, another to $\Tb$,
        the third to $\Tb'$ (belonging to an isotopic translate of
        $\Tb$), and following some fixed isotopy of $\Tb$ to $\Tb'$
        along the rounded edge. Stretching out the rounded edge, we
        obtain a chain homotopy between this map, and the map induced
        by counting pseudo-holomorphic triangles for the Heegaard
        triple $(\Sigma,\alphas,\betas,\betas',w,z)$, which is some
        cycle in $\CFKm(\Sigma,\betas,\betas',w,z)$ (which we can
        think of as the relative invariant of the isotopy). The fact
        that both maps induce Maslov grading-preserving isomorphisms
        on Floer homology ensures that the relative invariant
        represents the top--dimensional homology generator of
        $\HFKm(\Sigma,\betas,\betas',w,z)$. This completes the
        identification of the continuation map with the map induced by
        pseudo-holomorphic triangles, on the level of homology. (We
        have explained here the case where only the $\beta$ circles
        are moving; the case where both $\alpha$-- and
        $\beta$--circles are moving follows with straightforward
        notational changes.)
        
        Finally, observe that the map induced by the homeomorphism
        $F_t$ is also a closest-point map.  Thus, for all $t$
        sufficiently small, $A(F_t)$ acts by the identity on homology.
        
        By the compactness of $[0,1]$ we conclude that any $F$ which
        is isotopic to the identity acts by the identity map on knot
        Floer homology.  Consequently the action of $\DiffP(Y,K)$
        descends to an action of $\DiffP(Y,K)/\DiffP_0(Y,K)=\MCG(
        Y,K)$ on $\HFKm (Y,K)$, concluding the proof.
\end{proof}

Note that if $F\colon (Y_1,K_1,p_1)\longrightarrow (Y_2,K_2,p_2)$ is a
homeomorphism, then Definition~\ref{def:MapOfHomeo} can be adapted to
define a map $$\HFKm(Y_1,K_1,p_1)\longrightarrow \HFKm(Y_2,K_2,p_2),$$
which is well-defined up to the above action of the mapping class
group of $(Y_2,K_2)$.

\bibliographystyle{plain}
\bibliography{biblio}

\end{document}